\documentclass[preprint,11pt]{amsart}

\usepackage{amssymb}
\usepackage{amsthm}

\usepackage{amsmath}

\usepackage{graphicx}
\usepackage{subfigure}
\usepackage{bm}
\usepackage{color}
\newtheorem{remark}{Remark}
 \theoremstyle{definition}
\newtheorem{exmp}{Example}[section]

\DeclareMathOperator\supp{su pp}

\begin{document}

%\begin{frontmatter}



%\title{An efficient high-order Nystr\"{o}m scheme for acoustic scattering by inhomogeneous penetrable media}

\title[An efficient high-order Nystr\"{o}m scheme for inhomogeneous acoustic scattering]{An efficient high-order Nystr\"{o}m scheme for acoustic scattering by inhomogeneous penetrable media with discontinuous material interface}


\author{Akash Anand,  Ambuj Pandey,  B. V. Rathish Kumar, Jagabandhu Paul}

\maketitle

%\author[label1,label2]{Akash Anand}
%\ead{akasha@iitk.ac.in}
%\ead[url]{http://home.iitk.ac.in/~akasha}
%\address[label1]{Mathematics and Statistics, IIT Kanpur, Kanpur, UP 208016, India}
%\address[label2]{Computing and Mathematical Sciences, Caltech, Pasadena, CA 91125, USA}
%\author[label1]{Ambuj Pandey}
%\ead{ambuj@iitk.ac.in}
%\author[label1]{B.~V.~Rathish Kumar}
%\ead{bvrk@iitk.ac.in}
%\author[label1]{Jagabandhu Paul}
%\ead{jagapaul@iitk.ac.in}

\begin{abstract}
This text proposes a fast, rapidly convergent Nystr\"{o}m method for the solution of the 
Lippmann-Schwinger integral equation that mathematically models the scattering of time-harmonic
acoustic waves by inhomogeneous obstacles, while allowing the material properties to jump across the interface. The method works with overlapping coordinate charts
as a description of the given scatterer. In particular, it employs ``partitions of unity" to 
simplify the implementation of high-order quadratures along with suitable changes of parametric 
variables to analytically resolve the singularities present in the integral operator to achieve 
desired accuracies in approximations. To deal with the discontinuous material interface in a high-order manner, a specialized quadrature is used in the boundary region. The approach further utilizes an FFT based strategy that uses equivalent
source approximations to accelerate the evaluation of large number of interactions that arise in 
the approximation of the volumetric integral operator and thus achieves a reduced computational 
complexity of $O(N \log N)$ for an $N$-point discretization. A detailed discussion on the solution 
methodology along with a variety of
numerical experiments to exemplify its performance in terms of both speed and accuracy are presented in this paper.
\end{abstract}

%\begin{keyword}
%Acoustic scattering \sep Lippmann-Schwinger integral equation \sep high-order methods \sep integral equation methods
%%% keywords here, in the form: keyword \sep keyword
%
%%% PACS codes here, in the form: \PACS code \sep code
%\PACS 02.60.Cb \sep \PACS 02.70.-c \sep \PACS 43.20.+g \sep \PACS 44.05.+e
%
%%% MSC codes here, in the form: \MSC code \sep code
%%% or \MSC[2008] code \sep code (2000 is the default)
%
%\end{keyword}

%\end{frontmatter}


\section{Introduction}
\label{sec::intro}

The solution of direct acoustic scattering problem where the goal is to obtain the scattered wave produced 
as a result of an interaction between a given incident wave and a given bounded inhomogeneity continues to
constitute one of the most challenging problems in computational science, especially where applications of
interest require such computations to be carried out for large penetrable scatterers with complex geometries.
Formally, the $d$-dimensional direct acoustic scattering problem that we consider in this paper for $d = 2$ 
and $3$ is described as follows: given an obstacle $\Omega$, a bounded open subset of $\mathbb{R}^d$, with a smooth boundary $\partial\Omega$, and an 
incident time-harmonic acoustic wave $\mathfrak{u}^{i}$ satisfying
\begin{equation}
\Delta \mathfrak{u}^{i}(\bm{x}) + \kappa^2 \mathfrak{u}^{i}(\bm{x}) = 0, \ \ \bm{x} \in \mathbb{R}^d,
\end{equation}
where $\kappa = \omega / c_0$ is the wavenumber, $\omega$ is the angular frequency, and $c_{0}$ is the constant speed of wave outside the inhomogeneity $\Omega$,
find the total acoustic field  
$\mathfrak{u}$ that satisfies {\cite{colton2012inverse}}
\begin{equation} \label{HE}
\Delta \mathfrak{u}(\bm{x}) + \kappa^2 n^2(\bm{x}) \mathfrak{u}(\bm{x}) = 0, \ \  \bm{x} \in \mathbb{R}^d,
\end{equation}
with the refractive index $n(\bm{x}) = c_{0} / c(\bm{x})$, where $c$, the speed of acoustic wave, is allowed to vary with position within $\Omega$ 
and the scattered field $\mathfrak{u}^s := \mathfrak{u} - \mathfrak{u}^i$ satisfies Sommerfeld radiation condition %\cite [p-16,p-66]{colton2012inverse}
\begin{equation} \label{eq:-Sommerfeld}
\lim_{r \to \infty} r^{(d-1)/2} \left( \frac{\partial\mathfrak{u} ^s}{\partial r} - ik\mathfrak{u}^s \right) = 0,
\end{equation}
where $r=\|\bm{x}\|_{2} = \sqrt{\sum_{j=1}^d x_j^2}$.
%\medskip

One of the main challenges in obtaining numerical approximation of the solution to the scattering problem arise 
from the need to accurately describe highly oscillatory functions in large domains. 
As a certain fixed number of points is necessary to resolve a wavelength, we typically require a large computational 
grid to obtain any meaningful solution and any effort to further reduce the computational errors demands that the size 
of discretization grows proportionally in all dimensions. It is, therefore, highly desirable to devise numerical 
methodologies that are efficient as well as high-order accurate. In this direction, several methodologies,
based on diverse formulations of the problem, ranging from differential equation to  variational form and to 
integral equation, have been proposed by various authors. In the solution strategies that are based on differential
equation or their corresponding weak formulations \cite{bayliss1985accuracy,coyle2000scattering,gan1993finite,li2010coupling,
 kirsch1990convergence,kirsch1994analysis,meddahi2003computing,rachowicz2000hp,stupfel2001hybrid}, 
% condition (\ref{eq:-Sommerfeld}), 
a relatively large computational 
domain containing the scatterer must be used, together with appropriate absorbing boundary conditions on
the boundary of the computational domain. 
%It is also known that they suffer from the dispersion error \cite{bayliss1985accuracy} and hence the accuracy of the approximate solution deteriorates rapidly as wave number increases. The dispersion error could be reduced either by using a finer mesh, or by using higher order methods but both of these technique increases the overall computational cost.
Thus, these procedures often require a large number of unknowns and, hence, lead to large linear systems. In addition, accurate absorbing
boundary conditions with efficient numerical implementations are quite difficult to construct; the error
associated with such boundary conditions typically dominates the error in the computed solution. 

In contrast, the integral equation approach, where the mathematical formulation directly ensures that the solution satisfies condition (\ref{eq:-Sommerfeld}) 
by suitably employing the radiating fundamental solution, is free from considerations mentioned above, and consequently,
does not require solution strategies to discretize outside the inhomogeneity. We, therefore, base our numerical treatment of the scattering problem on an 
equivalent integral equation formulation which is given by the \textit{Lippmann-Schwinger} equation {\cite{colton2012inverse, martin2003acoustic},}
 \begin{equation} \label{eq:-Lippmann}
 \mathfrak{u}(\bm{x}) + \kappa^2 \int \limits_{\mathbb{R}^{d}}G_{\kappa}(\bm{x},\bm{y})m(\bm{y})\mathfrak{u}(\bm{y})d\bm{y} = 
 \mathfrak{u}^{i}(\bm{x}), \hspace{3mm} \hspace{3mm} \bm{x} \in \mathbb{R}^{d},
\end{equation}
where
\begin{equation} \label{eq:-G}
 G_{\kappa}(\bm{x},\bm{y}) =
\begin{cases}
\frac{i}{4}H^{1}_{0}(\kappa|\bm{x}-\bm{y}|), & \text{in } \mathbb{R}^{2}, \\
\exp({i\kappa|\bm{x}-\bm{y}| })/4\pi|\bm{x}-\bm{y}|, & \text{in }  \mathbb{R}^{3} ,
\end{cases}
\end{equation}
is the radiating fundamental solution of Helmholtz equation in the free space and $m(\bm{x}) = 1-n^{2}(\bm{x})$.

In recent years, a lot of progress has been made toward numerical solution of the Lippmann-"Schwinger equation; for example, see
\cite{aguilar2004high,anand2006efficient,anand2007efficient,andersson2005fast,chen2002fast,fan2001cgfft,hesford2010fast,hyde2002fast, liu2000applications,
	 polimeridis2014stable,vainikko2000fast,xu2001fast,zhu2003quadrature,zwamborn1992three, duan2009high, gillman2014spectrally}.
Most fast numerical schemes among these, though high order accurate for smooth scattering media,
exhibit only linear convergence in the presence of material discontinuity \cite{aguilar2004high,andersson2005fast,fan2001cgfft,vainikko2000fast, duan2009high, gillman2014spectrally}.
For example, in \cite{duan2009high}, Duan and Rokhlin, introduced an efficient high order quadrature formula which utilizes a corrected Trapezoidal
rule in conjunction with FFT for fast and accurate approximation of volume integration in (\ref{eq:-Lippmann}). 
The high order convergence of the scheme, however, requires the refractive index $n(\bm{x})$ to be globally smooth.  More recently, a high order direct solver with $O(N^{3/2})$ computational complexity has been proposed in \cite{gillman2014spectrally}. 
%The algorithm is hierarchical in nature, and split the problem into two, namely, a variable coefficient problem in an inhomogeneous medium and a constant coefficient problem in 
%homogeneous medium. In each of the domain solution is obtained by constructing Dirichlet-to-Neumann (DtN) map on the boundary of
%scatterer $\partial \Omega$ and final solution is obtained by merging these operator on the boundary $\partial \Omega$. 
Again, while this algorithm is robust and provide accurate approximation for smooth scatterer, it does not exhibit rapid convergence in the presence of a discontinuous material interface.
A couple of fast techniques, one that rely on the use of ``discontinuous FFT'' \cite{zhu2003quadrature} for efficient computations 
while the other uses FFT for accelerated evaluation of convolution in polar coordinates after a suitable decomposition of the 
Green's function via addition theorem \cite{hyde2002fast}, however, achieve second order accuracy for discontinuous scattering configurations.
Another approach for solution of acoustic volumetric scattering problem introduced
 in \cite{anand2006efficient,anand2007efficient}, though high-order convergent, is designed to be computationally efficient only for 
 ``thin'' scattering configurations. 
 This scheme gains high-order convergence through a combination of changes of parametric variables (in order to resolve the singularities present in the
 Green function) and by suitably employing ``partitions of unity" to yield smooth and periodic integrand away from vicinity of target points.
 Our present approach,
in fact, is a non trivial extension of ideas presented in %those references 
\cite{anand2006efficient,anand2007efficient} wherein
we obtain a solver for general scattering configurations that exhibits computational complexity of $O(N \log N)$ with respect to the grid size $N$ while 
retaining high-order accuracy even in the presence of material discontinuity.  We believe that this algorithm is the first fast integral equation solver which provide high-order convergence for full volumetric scattering problem with a discontinuous material interface.

The rest of the paper is organized as follows. In Section \ref{sec::algo}, we present main algorithmic components of our numerical scheme. We 
then present, in Section \ref{sec::2d}, a detailed account of our approximation strategy for the integral operator arising in the context of two 
dimensional volumetric scattering problem. The last part of this section presents a series of numerical experiments to exemplify the performance of
the proposed method, both in terms of accuracy and in terms of computational efficiency. This is followed,  in Section \ref{sec::3d},  by a brief
discussion on numerical solution of corresponding three dimensional scattering problems using a straightforward extension of our two dimensional
solver and demonstrate its effectiveness via some numerical experiments. Finally, our conclusions are summarized in Section \ref{sec::conclu}.

%In this paper, we present a fast iterative solver for \textit{Lippmann-Schwinger} equation 
%where acoustic waves are scattered by penetrable inhomogeneous media in two dimensions
%that can readily be extended to three dimensions. Notably our algorithm not only converges with high-order but
%also exhibit a computational complexity of $\mathcal{O}(N \log N)$. 

\section{Principal components of the method}
\label{sec::algo}

%This section presents key components of the proposed algorithm. 
%We discretize equation (\ref{eq:-Lippmann}) using Nystr\"{o}m scheme and the resulting linear system 
%is solved using GMRES. 

As mentioned earlier, we base our numerical strategy on solving the linear system arising out of a
Nystr\"{o}m 
discretization of the Lippmann-Schwinger integral equation. Given the denseness of resulting linear systems where use of direct linear 
solvers can be prohibitively expensive, we employ the matrix-free version of the iterative solver GMRES \cite{saad1986}, in fully complex arithmetic, for the solution of the discrete form of (\ref{eq:-Lippmann}). 
In this context, we focus our presentation on the computational technique for accurate and efficient evaluation of the integral operator 
\begin{equation} \label{eq:-VolInt}
 \mathcal{K}[\mathfrak{u}](\bm{x}) = \int \limits_{\Omega}G_{\kappa}(\bm{x},\bm{y})m(\bm{y})\mathfrak{u}(\bm{y})d\bm{y}.
\end{equation}
%directly results 
%We,
%therefore, focus our presentation 
%%of the numerical schemes 
%on an efficient integration scheme to
%approximate the integral (\ref{eq:-VolInt}) accurately. 

%Although integral equation method has a lot of advantages and it  provide the most accurate solutions to wave
%scattering problems but direct use of it would lead dense linear systems which require $\mathcal{O}(N^2)$ operation
%per iteration of an iterative solver--where $N$ is the number of unknowns. Therefore for large realistic problems 
%it become computationally prohibitive and still formidable even for super computers. 

\begin{figure}[h]
\begin{center}
\includegraphics[scale=0.7]{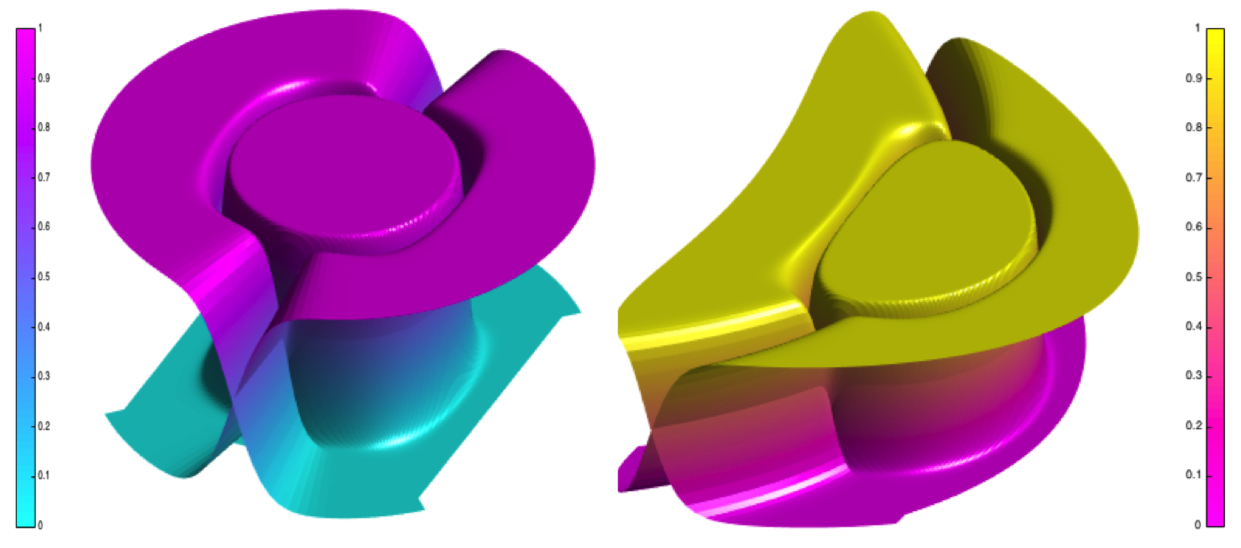}
%  \subfigure[POU on disc ]{\includegraphics[clip=true, trim=300 150  320 90, scale=0.20]{images/disc_Pou.png}} 
%  \hspace{10mm}
 % \subfigure[POU on bean]{\includegraphics[clip=true, trim=370 150  270 90, scale=0.20]{images/bean_Pou_7.png}}
%\subfigure[POU Over Bean]{\includegraphics[width=3.2in, height=2in,angle=0]{beanPou.eps}}
\caption{Partitions of unity.} \label{fig:-pou}
\label{pou}
\end{center}
\end{figure}

We start by describing the scatterer $\Omega$
through a collection $\mathcal{P} = \{ \mathcal{P}_{k}\}_{k=1}^{K}$ of $K$ overlapping coordinate patches
where the $k$-th patch is homeomorphic to an open set $\mathcal{H}_{k} \subset (0,1)^{d}$ via a smooth
invertible parametrization  $\bm{\xi}_{k}$.  This, in conjunction with a partitions of unity (POU) subordinated to the covering $\mathcal{P}$, that is, functions $\omega_{k}(\bm{x}):k=1,\ldots,K$, satisfying
\[
 \sum \limits_{k=1}^{K} \omega_{k}(\bm{x}) =1 \hspace{2mm} \text{\ for all \ } \bm{x} \in \Omega,
\]
where for each $k$, $\omega_{k} \in C^{\infty}(\Omega)$ with its support contained in $\mathcal{P}_{k}$,
reduces the evaluation of (\ref{eq:-VolInt})
to computation of integrals over these $K$ patches. % $\mathcal{P}_{\mathrm{k}}$.
Note here that for all those patches $\mathcal{P}_{k}$ whose
closure does not intersect with the boundary of $\Omega$, the corresponding $\omega_{k}$ vanishes to high order along with all of its
derivatives on the boundary of the patch. On remaining patches, however, where one of the edges coincide with the boundary of $\Omega$, $\omega_{k}$
clearly does not vanish and, in fact, attains the value $1$ (see Figure~\ref{fig:-pou} for examples).
For the clarity of presentation of the proposed quadrature, and to distinguish between patches based on this criterion,
we introduce two index sets, namely, $ \mathcal{I}_{I} = \{ k \ |\  \overline{\mathcal{P}_{k}} \cap \partial \Omega = \O{} \}$ and $ \mathcal{I}_{B} = \{k \ |\  \overline{\mathcal{P}_{k}} \cap \partial \Omega \neq \O{} \}$ 
corresponding to {\em interior} and {\em boundary} patches respectively.
%set of patches $\mathcal{P}$ into two sets, namely, set of boundary patches
%$\mathcal{P}_{B} = \Big\{\mathcal{P}_{\mathrm{k}}\in \mathcal{P}\ \big|\  \mathrm{k}\in \mathcal{I}_{B} \Big\}$
%and set of interior patches $\mathcal{P}^{I} = \Big\{\mathcal{P}_{\mathrm{k}} \in \mathcal{P}\big| \mathrm{k}\in \mathcal{I}_{I}\Big \}$
%where $ \mathcal{I}_{B} = \{\mathrm{k}| \mathcal{P}_{\mathrm{k}}\cap \partial \Omega \neq \O{} \}$ and
%$ \mathcal{I}_{I} = \{\mathrm{k}| \mathcal{P}_{\mathrm{k}}\cap \partial \Omega = \O{} \}$
%define an index set for boundary and interior patches respectively. Note that $ \mathcal{P}^{B} \cap  \mathcal{P}^{I} = \O{}$ and 
 %$ \mathcal{P}^{B} \cup  \mathcal{P}^{I} = \mathcal{P}$.  
Now, 
%using POU functions $\omega_{\mathrm{k}}$ 
clearly, the integral (\ref{eq:-VolInt}) can be reexpressed as a sum of integrals over boundary and interior patches
\begin{align} \label{eq:-VolInt-B-I}
 \mathcal{K}[\mathfrak{u}](\bm{x}) =  \sum\limits_{k \in \mathcal{I}_{B}} \int_{\mathcal{P}_{k}} G_{\kappa}(\bm{x},\bm{y})m(\bm{y})\mathfrak{u}(\bm{y})\omega_{k}(\bm{y})d\bm{y}
  +\sum\limits_{k \in \mathcal{I}_{I}} \int_{\mathcal{P}_{k}} G_{\kappa}(\bm{x},\bm{y})m(\bm{y})\mathfrak{u}(\bm{y})\omega_{k}(\bm{y})d\bm{y}.
\end{align}
%\begin{align} \label{eq:-VolInt-B-I}
% \mathcal{K}[\mathfrak{u}](\bm{x}) &=  \sum\limits_{\mathrm{k} \in \mathcal{I}_{B}} \int \limits_{\mathcal{P}_{\mathrm{k}}} G(\bm{x},\bm{y})m(\bm{y})\mathfrak{u}(\bm{y})\omega_{k}(\bm{y})d\bm{y} \nonumber \\
%  &+\sum\limits_{\mathrm{k} \in \mathcal{I}_{I}} \int \limits_{\mathcal{P}_{\mathrm{k}}} G(\bm{x},\bm{y})m(\bm{y})\mathfrak{u}(\bm{y})\omega_{\mathrm{k}}(\bm{y})d\bm{y}.
%\end{align}
%Finally, the goal is to evaluate each integral in the above equation accurately and efficiently
%at each Nystr\"{o}m node points. 
In view of this, we introduce a set of quadrature points on the $k$-th patch, say $\mathfrak{T}_k$, by
taking the image of the equi-spaced  grid $\{(i_1/N_1,\ldots, i_{d-1}/N_{d-1},i_d/N_d) \ | \ 0 \le i_{\ell} \le N_{\ell}, \ell = 1,\ldots,d \}$ in $[0,1]^d$ %(where the actual integration is performed), 
under the map $\bm{\xi}_{\mathrm{k}}$. 
Further, the union of these grid points, that is,
\begin{equation}
 \mathfrak{T} = \bigcup \limits_{k=1}^{K}
\mathfrak{T}_k
\end{equation}
where
\begin{equation}
\mathfrak{T}_k = 
\left\{ \bm{\xi}_{k}(i_1/N_1,\ldots, i_{d-1}/N_{d-1},i_d/N_d) \ | \ 0 \le i_{\ell} \le N_{\ell}, \ell = 1,\ldots,d \right\} 
\end{equation}
defines the complete set of Nystr\"{o}m discretization points where we seek to approximate the solution of equation (\ref{eq:-Lippmann}).

Indeed, the performance %in employing an
of our iterative solver hinges on our ability to
approximate the integrals in (\ref{eq:-VolInt-B-I}) accurately, in a computationally efficient manner. In the next section, %(\ref{sec::2d})
we provide a detailed description of a fast, high order integration scheme %to evaluate integrals in equation (\ref{eq:-VolInt-B-I}) 
in two dimensions
that can be readily extended to three dimensions, which we briefly outline in Section \ref{sec::3d}.
\begin{remark}
In this text, we often refer to $\bm{x}$ in equation (\ref{eq:-VolInt-B-I}) as a {\em target point} whereas points $\bm{y}$ therein have sometimes been called {\em source points}. 
%\item For boundary patches, we identify the $d$-th coordinate variable in $[0,1]^d$ with the transverse parameter. We, therefore, use the notation $(\bm{u},t) =(u_{1},\ldots u_{d-1},t)$ to denote a point in the corresponding parameter space and assume that the boundary coincides with $t = 1$.
%\item $(u_{0},t_{0}) =(u_{0,1},\cdots u_{0,d-1},t_{0}) $ throughout in this paper to denote the coordinate of target point $\bm{x}$
% in the parametric space $\mathcal{H}_{\mathrm{k}}\subset \mathbb{R}^{d}$.
% \item In the two dimension we take $\mathcal{H}_{\mathrm{k}} = (0,1)^2$ and in three dimension
% $\mathcal{H}_{\mathrm{k}} = (0,1)^3$.
\end{remark}
\begin{remark}
While discussing the approximations on boundary patches, we identify the $d$-th coordinate variable in $[0,1]^d$ with the transverse parameter. We,
therefore, use the notation $\bm{t}=(t_{1},\ldots, t_{d-1},t_d)$ to denote a point in the corresponding parameter space and assume that the boundary coincides with $t_d = 0$. 
%We also note for future reference that the parametric coordinate of target point $\bm{x}$ are denoted by $\bm{u}_0$ which, in the context of boundary parches read $(\bm{u}_{0},t_{0}) =(u_{0,1},\cdots u_{0,d-1},t_{0}) $. 
\end{remark}

\section{Fast and accurate evaluation of integrals in two dimensions} 
\label{sec::2d}

The difficulty in accurately computing integrals in (\ref{eq:-VolInt-B-I}) is significantly more when the target point $\bm{x}$ 
lies in the integration patch $\mathcal{P}_k$ compared to the case when it does not. Indeed, when $\bm{x} \in \mathcal{P}_k$, 
owing to the singularity of the kernel $G_{\kappa}(\bm{x},\bm{y})$ at $\bm{y} = \bm{x}$, the integrand is unbounded within the integration 
domain and direct use of a standard quadratures yield inaccurate approximations. Thus, specialized quadrature rules must be
developed and used to deal with such singular integrals. The case when $\bm{x} \not\in \mathcal{P}_k$, in contrast, does not present this challenge.
To effectively present these two contrasting scenarios, we further refine the index sets $\mathcal{I}_{B}$ and $\mathcal{I}_{I}$ 
by introducing target point dependent index sets $\mathcal{M}_B(\bm{x}) = \left \{ k \in \mathcal{I}_{B} \ |\ \bm{x} \in \mathcal{P}_k  \right\}$ and $\mathcal{M}_I(\bm{x}) = \left \{k \in \mathcal{I}_{I} \ |\ \bm{x} \in \mathcal{P}_{k} \right\}$ to 
rewrite (\ref{eq:-VolInt-B-I}) as
\begin{align}\label{eq:-VolInt-B-II}
 \mathcal{K}[\mathfrak{u}](\bm{x}) &= \left( \sum\limits_{k \in \mathcal{M}_B(\bm{x})} + \sum\limits_{k \not\in \mathcal{M}_B(\bm{x})} \right) \int _{\mathcal{P}_{k}}
 G_{\kappa}(\bm{x},\bm{y}) m(\bm{y})\mathfrak{u}(\bm{y})\omega_{k}(\bm{y})\,d\bm{y} \nonumber \\ &+ \left( \sum\limits_{k \in \mathcal{M}_I(\bm{x})} + \sum\limits_{k \not\in \mathcal{M}_I(\bm{x})} \right) \int _{\mathcal{P}_{k}}
 G_{\kappa}(\bm{x},\bm{y}) m(\bm{y})\mathfrak{u}(\bm{y})\omega_{k}(\bm{y})\,d\bm{y}.
\end{align}
Clearly, each of the integrals above can be rewritten in the parametric coordinates
\[
\int \limits_{\mathcal{P}_{k}}
 G_{\kappa}(\bm{x},\bm{y}) m(\bm{y})\mathfrak{u}(\bm{y})\omega_{k}(\bm{y})\,d\bm{y} = 
 \iint \limits_{[0,1]^2}
 G_{\kappa}(\bm{x},\bm{\xi}_k(t_1,t_2)) \varphi_k[\mathfrak{u}](t_1,t_2) \xi'_{k}(t_1,t_2) \,dt_1 \, dt_2,
\]
where
$
\varphi_k[\mathfrak{u}](t_1,t_2) = m(\bm{\xi}_{k}(t_1,t_2))\mathfrak{u}(\bm{\xi}_{k}(t_{1},t_{2}))\omega_{k}(\bm{\xi}_{k}(t_1,t_2)),
$
and $\xi'_{k}$ denotes the Jacobian of the transformation $\bm{\xi}_{k}$. 

Now, for the cases $k \not\in \mathcal{M}_B(\bm{x})$ and $k \not\in \mathcal{M}_I(\bm{x})$ in (\ref{eq:-VolInt-B-II}), %the target point $neighborhood{x}$ do not lie in the integration patch which in turn ensures that 
$G_{\kappa}(\bm{x},\bm{y})$ remains non-singular throughout the region of integration. Though adopting a single high-order approximation
scheme for both scenarios is possible, we actually utilize two different quadratures to take advantage of a more favorable behavior of the integrands in the 
later case when $k \not\in \mathcal{M}_I(\bm{x})$ where $\varphi_k[\mathfrak{u}]$ vanish to high order at the boundary of the integration domain and have smooth periodic extensions to $\mathbb{R}^2$. 
%In further discussions, we call these two generic types of integrals as {\bf RS-TYPE} and {\bf RP-TYPE} integrals respectively, indicating a {\em Regular-Semiperiodic} integrand for the case $k \in \mathcal{M}_B(\bm{x})$ and {\em Regular-Periodic} integrand in the case $k \in \mathcal{M}_I(\bm{x})$.
As is well known, the trapezoidal rule exhibits super-algebraic convergence for 
smooth and periodic integrands, which we indeed employ to obtain accurate approximations in this case.
In contrast, when $k \not\in \mathcal{M}_B(\bm{x})$, a straightforward use of trapezoidal rule does not produce high-order 
accuracy as the integrands do not vanish at $t_2 = 0$. We overcome this minor difficulty by utilizing the trapezoidal rule for 
the integration with respect to $t_1$-variable while employing a composite Newton-Cotes quadrature in transverse variable $t_2$ to achieve 
approximations whose rate of convergence directly depends on the order of Newton-Cotes used and could be enhanced arbitrarily as long as the smoothness of $m$ within $\Omega$ allows it. 
%However, on the other hand, integrands
% of the integrals in the fourth term (integral over interior patches) are smooth and vanishes to high-order
%at the end points of integration interval in each variable therefore high-order can be obtained by means of 
%repeated use of Trapezoidal rule quadrature. Similar to the previous case, this case also requires only
% $\mathcal{O}(\frac{N^{2}}{L^{2}})$ operations.

As mentioned above, integrands in (\ref{eq:-VolInt-B-II}) corresponding to the cases when $k \in \mathcal{M}_B(\bm{x})$ and $k \in \mathcal{M}_I(\bm{x})$ are singular. 
%that we approximate using specialized quadratures. 
To achieve rapidly convergent approximations, we rely on analytic resolution of singularities through suitable changes of parametric 
variables and application of high-order quadratures to resulting smooth integrands. 
%We discuss the details of this approach in section \ref{sec:-SI2D}.
In addition, we also localize the region where such coordinate transformations are affected to a small neighborhood of the singular
point using a suitable smooth and compactly supported cut-off function. Indeed, we see that
\begin{equation} \label{eq:-VolInt-B-III}
 \iint \limits_{[0,1]^2}
 G_{\kappa}(\bm{x},\bm{\xi}_k(\bm{t})) \cdots \,d\bm{t} = 
  \iint \limits_{[0,1]^2}
 G_{\kappa}(\bm{x},\bm{\xi}_k(\bm{t})) \cdots \eta(\bm{t};\bm{\xi}_k^{-1}(\bm{x})) \,d\bm{t} +
  \iint \limits_{[0,1]^2}
 G_{\kappa}(\bm{x},\bm{\xi}_k(\bm{t})) \cdots (1-\eta(\bm{t};\bm{\xi}_k^{-1}(\bm{x}))) \,d\bm{t},
\end{equation}
where $\eta(\cdot;\bm{t}_0)$ is a $C^{\infty}$ function that it is compactly supported in a neighborhood of 
the point $\bm{t}_0 \in [0,1]^2$ while $\eta \equiv 1$ in a smaller neighborhood of $\bm{t}_0$. This localization of 
singularity as seen in the first integral on the right hand side of (\ref{eq:-VolInt-B-III}) brings in a two-fold benefit,
namely, {\em (a)} it limits the relatively expensive treatment of singularity which, among other computational challenges,
also demands an interpolation of grid-data to off-grid quadrature points to a small integration domain, and {\em (b)}
it allows for additional speed-up in the computation of non-singular second term on the right hand side of (\ref{eq:-VolInt-B-III}).
Note that this regular integral, of course, can be integrated to high-order using the same numerical quadratures that
we apply in the case of $\bm{x} \not\in \mathcal{P}_k$, as explained above.
%Again, instead of using a single quadrature for to deal with both, we utilize two different approximation schemes, especially, to take advantage of the vanishing nature of the integrand in the case of interior patches. In further descriptions, we call these two generic types of integrals as {\bf SS-Type} and {\bf SP-Type} integrals, respectively, indicating a {\em Singular-Semiperiodic} integrand for the case $k \in \mathcal{M}_B(\bm{x})$ and {\em Singular-Periodic} integrand in the case $k \in \mathcal{M}_I(\bm{x})$. We present a detailed discussion on accurate approximation of these singular integrals in section \ref{sec:-SI2D}.

It is straightforward to see that the direct application of the integration scheme that we described above leads to a computational complexity of
$O(N^2)$, where $N = | \mathfrak{T}|$ is the size of the set $ \mathfrak{T}$, the total number of quadrature points in the Nystr\"{o}m scheme. However, 
as we outline below, one can improve the computational complexity of this methodology to $O(N \log N)$ by breaking the overall computation of these integrals into two parts -- a relatively
small but specialized calculation in the neighborhood of singularity and remaining non-singular contributions arising from the voluminous bulk.
Toward this, we begin by bounding the inhomogeneity $\Omega$ by a square cell $\mathcal{C}$ of side length $A$.
%The square cell $\mathcal{C}$ 
This cell is further partitioned into $L^{2}$ identical cells $\mathcal{C}_{ij}$ $(i,j =1,\cdots,L )$ of side length $H = A/L$ such that
the bounding cell $\mathcal C$ contains $L$ of these smaller cells along its sides (for an example, see Figure \ref{fig:-accel}). %and the subscript $ij$ identifies each of these cells $c_{ij}$ uniquely.
%Thus, on an average, each cell contain $\mathcal{O}(\frac{N}{L^{2}})$ discretization points. 
We assume, without loss of generality, that for the chosen side length $A$, cells $\mathcal{C}_{ij}$ do not admit
inner acoustical resonance. This, of course, can be easily ensured by necessary adjustment in the choice of $A$ so that $-{\kappa}^{2}$ 
is not a Dirichlet eigenvalue of the Laplace operator in the cell $\mathcal{C}_{ij}$.
%Now we introduce concept 
These cells are used to break the computation of integrals appearing in equation (\ref{eq:-VolInt-B-II}) into adjacent and non-adjacent 
calculations which, as pointed out above, is primarily motivated by our desire for the method to have a more favorable computational 
complexity than the quadratic cost in the number of unknowns.
%This idea of adjacency can be summarized in the following manner: 

We say a source point $\bm{y}$ is {\em adjacent} to a fixed point target point $\bm{x} \in \mathcal{C}_{ij}$, if it belongs to the set $\mathcal{N}(\bm{x})$ defined by
\begin{equation}\label{setdef:-Adj}
 \mathcal{N}(\bm{x}) = \left \{\bm{y}\ |\ \bm{y} \in C_{kl} \text{ for some } k,l \text{ satisfying }  | k-i | \leq 1, | l-j | \leq 1  \right \}.
\end{equation}
Obviously, we say $\bm{y}$ is {\em non-adjacent} to $\bm{x}$ if $\bm{y} \not\in \mathcal{N}(\bm{x})$.
%Then clearly $\mathcal{N}_{x}$ constitute a square region of side length $3H$ containing cell $c_{i_{0},j_{0}}$ and its $8$ neighboring cells.
%Now any point $\bm{y} \in c_{ij}$ are called adjacent to the target point $\bm{x} \in c_{i_{0},j_{0}}$ if $\bm{y} \in \mathcal{N}_{x}$,
%and nonadjacent if otherwise. In this context we refer adjacent set $\mathcal{N}_{x}$ as a super cell.
Based on this, we separate the integrals over $\mathcal{P}_k$ in (\ref{eq:-VolInt-B-II}) into integrals over 
$\mathcal{P}_k \cap \mathcal{N}(\bm{x})$ accounting for the adjacent contributions and the computationally 
large non-adjacent contributions for which a variant of the two face FFT based acceleration strategy 
introduced in \cite{bruno2001fast} is employed for efficient calculations. 
%The difficulty in accurately computing the adjacent contribution significantly increases when the target point $\bm{x}$ lies in the integration patch $\mathcal{P}_k$ making the integrand singular compared to the case when it does not. To effectively present these two contrasting scenarios, we further refine the index sets $\mathcal{I}_{B}$ and $\mathcal{I}_{I}$ by introducing target point dependent index sets $\mathcal{M}_B(\bm{x}) = \left \{ k \in \mathcal{I}_{B} \ |\ \bm{x} \in \mathcal{P}_k  \right\}$ and $\mathcal{M}_I(\bm{x}) = \left \{k \in \mathcal{I}_{I} \ |\ \bm{x} \in \mathcal{P}_{k} \right\}$. 
%
%Thus, rewriting (\ref{eq:-VolInt-B-I}) as
%\begin{align}\label{eq:-VolInt-B-II}
% \mathcal{K}[\mathfrak{u}](\bm{x}) &= \left( \sum\limits_{k \in \mathcal{M}_B(\bm{x})} + \sum\limits_{k \not\in \mathcal{M}_B(\bm{x})} \right) \int _{\mathcal{P}_{k}}
% G_{\kappa}(\bm{x},\bm{y}) m(\bm{y})\mathfrak{u}(\bm{y})\omega_{k}(\bm{y})\,d\bm{y},
%\end{align}
%
%Note that, owing to the singularity of the kernel $G(\bm{x},\bm{y})$ at $\bm{y} = \bm{x}$, the integrands in (\ref{eq:-VolInt-B-II}) corresponding to $k \in \mathcal{M}_B(\bm{x})$ and $k \in \mathcal{M}_I(\bm{x})$ are singular as, in both cases, the target point $\bm{x}$ belongs to the integration patch $\mathcal{P}_k$. 

We, thus, split integrals in (\ref{eq:-VolInt-B-II}) to rewrite the expression for $\mathcal{K}[\mathfrak{u}](\bm{x})$ as
\begin{align}\label{eq:-VolInt-Adj-NonAdj}
 \mathcal{K}[\mathfrak{u}](\bm{x}) &= \left\{ \left( \sum\limits_{k \in \mathcal{M}_B(\bm{x})} + \sum\limits_{k \not\in \mathcal{M}_B(\bm{x})} \right) \int _{\mathcal{P}_{k} \cap\mathcal{N}(\bm{x})} + \sum\limits_{k \in \mathcal{I}_{B}} \int _{\mathcal{P}_{k} \cap\mathcal{N}(\bm{x})^c} \right\}
 G_{\kappa}(\bm{x},\bm{y}) m(\bm{y})\mathfrak{u}(\bm{y})\omega_{k}(\bm{y})\,d\bm{y} \nonumber \\
 &+ \left\{ \left( \sum\limits_{k \in \mathcal{M}_I(\bm{x})} + \sum\limits_{k \not\in \mathcal{M}_I(\bm{x})} \right) \int _{\mathcal{P}_{k} \cap\mathcal{N}(\bm{x})} + \sum\limits_{k \in \mathcal{I}_{I}} \int _{\mathcal{P}_{k} \cap\mathcal{N}(\bm{x})^c} \right\}
 G_{\kappa}(\bm{x},\bm{y}) m(\bm{y})\mathfrak{u}(\bm{y})\omega_{k}(\bm{y})\,d\bm{y},
\end{align}
where $\mathcal{N}(\bm{x})^{c}$ denotes the compliment of $\mathcal{N}(\bm{x})$, that is, the set $\mathfrak{T} \setminus \mathcal{N}(\bm{x})$.

We note here that the computational cost of carrying out the adjacent calculations for all target points $\bm{x} \in  \mathfrak{T}$,
in view of the fact that on an average each cell contain $O\left(N/L^{2}\right)$ discretization
points, is  $N \times O\left(N/L^{2}\right) = O\left(N^{2}/L^{2}\right)$.

The cost involved in computing the contributions to integrals coming from $\mathcal{P}_k \cap \mathcal{N}(\bm{x})^{c}$ for all target points $\bm{x} \in \mathfrak{T}$, as we explain later in 
Section (\ref{subsec:-Non-Sing}) is given by 
$
O(LN^{1/2}\log (LN^{1/2}))+O(N^{3/2}/L^3)+O(N^{3/2}/L).
$
Thus, we get
$$O\left(LN^{\frac{1}{2}}\log \left(LN^{\frac{1}{2}}\right)\right)+O\left(\frac{N^{\frac{3}{2}}}{L^3}\right)+O\left(\frac{N^{\frac{3}{2}}}{L}\right)+O\left(\frac{N^2}{L^2}\right)$$ as
the total cost of computing $\mathcal{K}[\mathfrak{u}]$ at all Nystr\"{o}m nodes. This, of course, suggests that 
by choosing the parameter $L= O(N^{1/2})$, the computational complexity of the algorithm reduces to the desired $O(N\log N)$.

\subsection{Singular Integration} \label{sec:-SI2D}
Recall that when the target point $\bm{x}$ belongs to the integration patch $\mathcal{P}_k$, we break the integral as follows:
\begin{align} \label{eq:-VolInt-Adj-NonAdj-I}
&  \int \limits_{\mathcal{P}_{k}}G_{\kappa}(\bm{x},\bm{y}) m(\bm{y})\mathfrak{u}(\bm{y})\omega_{k} \,d\bm{y}  = \iint \limits_{[0,1]^2}
 G_{\kappa}(\bm{x},\bm{\xi}_k(\bm{t})) \varphi_k[\mathfrak{u}](\bm{t}) \xi'_k(\bm{t})(\bm{y}) \eta(\bm{t};\bm{\xi}_k^{-1}(\bm{x})) \,d\bm{t} 
 +   \nonumber \\  
&  \iint \limits_{\bm{\xi}^{-1}_k(\mathcal{P}_{k} \cap\mathcal{N}(\bm{x}))}
 G_{\kappa}(\bm{x},\bm{\xi}_k(\bm{t})) \varphi_k[\mathfrak{u}](\bm{t}) \xi'_k(\bm{t})(\bm{y}) (1-\eta(\bm{t};\bm{\xi}_k^{-1}(\bm{x}))) \,d\bm{t} +  \nonumber \\ 
 & \int \limits_{\mathcal{P}_{k} \cap\mathcal{N}(\bm{x})^c} G_{\kappa}(\bm{x},\bm{y}) m(\bm{y})\mathfrak{u}(\bm{y})\omega_{k}(\bm{y}) (1-\eta(\bm{\xi}_k^{-1}(\bm{y});\bm{\xi}_k^{-1}(\bm{x}))) \,d\bm{y}.
\end{align}

As mentioned above, we rely on the {\em FFT-accelerator} for computation of the third term in (\ref{eq:-VolInt-Adj-NonAdj-I}). This, however, entails
that we compute this integral as a convolution which, in turn, requires that the integrand remains independent of localization cut-off function $\eta$. 
This, of course, can easily be achieved by ensuring that the support of $\eta(\cdot,\xi_k^{-1}(\bm{x}))$ is contained in $\mathcal{N}(\bm{x})$. The second term in this expression, which we 
refer to as {\em local correction} to the accelerated numerics, is computationally small. We delay a more detailed discussion on these aspects to Section (\ref{subsec:-Non-Sing}). 
We detail the approximation of the singular integration
\begin{equation} \label{eq:-VolInt-I-P}
 \mathcal{K}_{k}^{sing}[\mathfrak{u}](\bm{x}) =   \iint \limits_{[0,1]^2\ \cap\ \supp\eta}
 G_{\kappa}(\bm{x},\bm{\xi}_k(\bm{t})) \varphi_k[\mathfrak{u}](\bm{t}) \xi'_k(\bm{t})(\bm{y}) \eta(\bm{t};\bm{\xi}_k^{-1}(\bm{x})) \,d\bm{t},
 \end{equation}
the first term on the right hand side of equation (\ref{eq:-VolInt-Adj-NonAdj-I}), next.

%Due to varying nature of integrands on the side of boundary and interior
%patches, we have used different strategy for evaluation of singular integral on boundary and interior patches--
%that we describe in what follows.

\subsubsection{Singular integration over interior patches} \label{sec:-Sing-Interior}
In this section, we describe high-order quadrature rule for evaluation of 
$ \mathcal{K}_{k}^{sing}[\mathfrak{u}](\bm{x})$ when $k \in \mathcal{M}_{I}(\bm{x}$). 
%In this scenario the
%integration patch $\mathcal{P}_{\mathrm{k}}$ happens to be an interior patch, i.e., $\mathrm{k} \in \mathcal{M}_{x}^{I}$ 
%and source point $\bm{y}\in \mathcal{P}_{\mathrm{k}} \cap\mathcal{N}_{x}$. 
%In order to simplify  our presentation we define for $ \mathrm{k}\in \mathcal{M}_{x}^{I}$ 
%We introduce the symbol $\mathcal{K}_k^I$ to denote the singular term
%\begin{equation} \label{eq:-VolInt-I-P}
% \mathcal{K}_{k}^{I}[\mathfrak{u}](\bm{x}) =   \iint \limits_{[0,1]^2\ \cap\ \supp\eta}
% G_{\kappa}(\bm{x},\bm{\xi}_k(\bm{t})) \varphi_k[\mathfrak{u}](\bm{t}) \xi'_k(\bm{t})(\bm{y}) \eta(\bm{t};\bm{\xi}_k^{-1}(\bm{x})) \,d\bm{t}.
%\end{equation}
%Similar to previous subsection again we start by localizing the kernel singularity by means of smooth cut-off function $\eta^{I}$
%as follows:
%\begin{align} \label{eq:-VolInt-Sing-Int}
% \mathcal{K}_{\mathrm{k}}^{I}[\mathfrak{u}](\bm{x}) &=  \int \limits_{0}^{1} \int \limits_{0}^{1} G(\bm{x},\bm{y}(u,t))\eta_{\bm{x}(u_{0},t_{0})}^{I}(u,t) \varphi_{\mathrm{k}}(u,t) du dt \nonumber \\
%   &+\int \limits_{0}^{1} \int \limits_{0}^{1} G(\bm{x},\bm{y}(u,t))\left(1-\eta_{\bm{x}(u_{0},t_{0})}^{I}(u,t)\right) \varphi_{\mathrm{k}}(u,t) du dt,
%\end{align}
%where function $\eta_{\bm{x}(u_{0},t_{0})}^{I}$ is defined as $
% \eta_{\bm{x}(u_{0},t_{0})}^{I} = \chi\left( |(u,t)-(u_{0},t_{0})| \right)$ and $|.|$ denotes euclidean norm in $\mathbb{R}^{2}$.
In this case, following \cite{bruno2001fast}, we use the localization function $\eta(\bm{t};\bm{\xi}_k^{-1}(\bm{x})) = \chi(|\bm{t}-\bm{\xi}_k^{-1}(\bm{x})|/r)$ with 
an appropriate real number $r > 0$ and a smooth function $\chi$ such that $\chi \equiv 1$ in a neighborhood of $0$ and $\chi(s) = 0$ for all $ s \ge 1$. 
%As we noted above, $\varphi_{k}[\mathfrak{u}]$ vanishes to high order
%in all directions on the boundary of interior patches $\mathcal{P}_{\mathrm{k}}$, $ \mathrm{k}\in \mathcal{I}_{I}$, hence the whole integrand is periodic in its domain of integration.   gvjyds
%Since the function $(1-\eta^{I})$ vanishes to high order in the vicinity of target point $\bm{x}$ therefore the
%whole integrand in the second integral of the equation (\ref{eq:-VolInt-Sing-Int}) can be viewed as a smooth and periodic functions
%which can be integrated to high order by means of repeated use of trapezoidal rule quadrature.
%However, the first integral, on the other hand, although enjoys a reduced domain of integration but contains weakly singular kernel.
We extend the domain of integration in (\ref{eq:-VolInt-I-P}) to the disc of radius $r$ around $\bm{\xi}_k^{-1}(\bm{x}) = (t_1^{\bm{x}},t_2^{\bm{x}})$ for each 
target point $\bm{x} \in \mathcal{P}_k \cap \mathfrak{T}$. Obviously, because of the fact that $\varphi_{k}[\mathfrak{u}]$ vanishes to high order in all directions 
on the boundary of the integration domain, this process does not compromise on the smoothness of the integrand.  To perform the integration,
we then change to polar coordinates centered at $(t_1^{\bm{x}},t_2^{\bm{x}})$:
\begin{equation*}
t_1 = t_1^{\bm{x}}+\rho \cos \theta, \hspace{2mm} t_2 = t_2^{\bm{x}} + \rho \sin \theta.
\end{equation*}
If we let $\tilde{\varphi}_{k}[\mathfrak{u}](t_1^{\bm{x}}+\rho \cos \theta,t_2^{\bm{x}}+\rho \sin \theta) = $
\begin{equation*}
 \varphi_{\mathrm{k}}[\mathfrak{u}](t_1^{\bm{x}}+\rho \cos \theta,t_2^{\bm{x}}+\rho \sin \theta)\xi'_k(t_1^{\bm{x}}+\rho \cos \theta,t_2^{\bm{x}}+\rho \sin \theta) \eta((t_1^{\bm{x}}+\rho \cos \theta,t_2^{\bm{x}}+\rho \sin \theta);(t_1^{\bm{x}},t_2^{\bm{x}})),
\end{equation*}
then the integral (\ref{eq:-VolInt-I-P}) takes the form
\begin{equation} \label{eq:-K-Int-Sing-Polar}
 \mathcal{K}_{\mathrm{k}}^{sing}[\mathfrak{u}](\bm{x}) = \frac{1}{2}\int \limits_{0}^{2\pi} \,d\theta \int  \limits_{-r_{1}}^{r_{1}}  |\rho|  G_{\kappa}(\bm{x},\bm{\xi}_k(t_1^{\bm{x}}+\rho \cos \theta,t_2^{\bm{x}}+ \rho \sin \theta))
\tilde{\varphi}_{k}[\mathfrak{u}](t_1^{\bm{x}}+\rho \cos \theta,t_2^{\bm{x}}+\rho \sin \theta) \,d\rho.
\end{equation}
The appearance of additional factor $|\rho|$ (the Jacobian of the polar change of variables) in the
integrand of integral (\ref{eq:-K-Int-Sing-Polar}) cancels the kernel singularity 
as $ |\rho| G_{\kappa}(\bm{x},\bm{y}(t_1^{\bm{x}}+\rho \cos \theta,t_2^{\bm{x}}+ \rho \sin \theta))$, clearly, is a smooth function of $\rho$. Additionally, the cut-off function $\eta$ vanishes to high order at the boundary of integration interval in $\rho$ variable, and therefore, a use of Trapezoidal rule 
 yields super-algebraic convergence for the $\rho$-integration. Moreover, the $\theta$-integral can also be approximated to high order by employing Trapezoidal rule quadrature as the corresponding integrand varies smoothly and periodically with respect to $\theta$.
 
While the change to polar coordinates provides a way to resolve the singularity analytically, the proposed application of
Trapezoidal rule demands that we provide the values $\tilde{\varphi}_{k}[\mathfrak{u}]$ at points outside of the computational grid $\mathfrak{T}$. 
This necessitates employing an efficient and accurate
interpolation strategy for evaluation of $\tilde{\varphi}_{k}[\mathrm{u}]$ at these newly transformed grid points.
Toward this, we adapt the Fourier refined polynomial interpolation introduced in \cite{bruno2001fast}, with a suitable choice of {\em polynomial degree},
to retain high order accuracy while maintaining computational efficiency.

\subsubsection{Singular integration over boundary patches} \label{sec:-Sing-Boundary}

When target point $\bm{x}$ lies in one of the boundary integration patches, say $\mathcal{P}_{k}$ 
with $k \in \mathcal{M}_{B}(\bm{x})$, we adopt a different approximation strategy to evaluate
$ \mathcal{K}_{k}^{sing}[\mathfrak{u}](\bm{x})$.
%
%\begin{equation} \label{eq:-VolInt-B-P}
% \mathcal{K}_{\mathrm{k}}^{B}[\mathfrak{u}](\bm{x}) =   \iint \limits_{[0,1]^2\ \cap\ \supp\eta}
% G_{\kappa}(\bm{x},\bm{\xi}_k(\bm{t})) \varphi_k[\mathfrak{u}](\bm{t}) \xi'_k(\bm{t})(\bm{y}) \eta(\bm{t};\bm{\xi}_k^{-1}(\bm{x})) \,d\bm{t}.
%\end{equation}
The first divergence comes in the form of the choice of the cut-off function $\eta$ where a circular support used for interior patches becomes unsuitable for use near the physical boundary of the domain, that is, near $t_2 = 0$. In this case, therefore, following the ideas introduce in \cite{anand2006efficient}, we take the localization function to have a rectangular support around $\bm{\xi}_k^{-1}(\bm{x}) = (t_1^{\bm{x}},t_2^{\bm{x}})$, given by
$$\eta((t_1,t_2);(t_1^{\bm{x}},t_2^{\bm{x}})) = \chi(| t_1- t_1^{\bm{x}} |/r)\chi(| t_2- t_2^{\bm{x}} |/r).$$
%We again decompose integral operator (\ref{eq:-VolInt-B-P}) in the form:
%\begin{align*}
%\mathcal{K}_{\mathrm{k}}^{B}[\mathfrak{u}](\bm{x}) &=\int \limits_{0}^{1} \int \limits_{0}^{1} G(\bm{x},\bm{y}(u,t))\eta_{\bm{x}(u_{0},t_{0})}^{B}(u,t) \varphi_{\mathrm{k}}(u,t) du dt \\
%&+\int \limits_{0}^{1} \int \limits_{0}^{1} G(\bm{x},\bm{y}(u,t))\left(1-\eta_{\bm{x}(u_{0},t_{0})}^{B}(u,t)\right) \varphi_{\mathrm{k}}(u,t) du dt. 
%\end{align*}
%The Integrand in the second integral of above equation is smooth, since $(1-\eta_{\bm{x}(u_{0},t_{0})})$ vanishes to 
%high-order in the vicinity of $(u_{0},t_{0})$ and it enjoys reduced region of integration because source point
%$\bm{y}\in \mathcal{P}_{k} \cap\mathcal{N}_{x}$. Therefore, it can be integrated efficiently and accurately
%using Trapezoidal rule quadrature in tangential direction (in $u$) and high order composite newton cotes quadrature in 
%transverse direction (in $t$). However, at the
%same time integrand in the first integral is still singular. In order to treat this singularity we follow 
%the ideas introduced in \cite{anand2006efficient}--- that we describe in what follows. 
We represent $ \mathcal{K}_{k}^{sing}[\mathfrak{u}](\bm{x})$ as
 \begin{equation} \label{eq:-VolInt-Adj-Sing-1}
\mathcal{K}_{k}^{sing}[\mathfrak{u}](\bm{x}) = \int \limits_{0}^{1} \mathcal{J}_{k}(t_2;\bm{x}) \chi(| t_2- t_2^{\bm{x}} |/r) dt_2
\end{equation}
where
\begin{equation} \label{eq:-I}
 \mathcal{J}_{k}(t;\bm{x}) = \int \limits_{t_1^{\bm{x}}-r}^{t_1^{\bm{x}}+r}G_{\kappa}(\bm{x},\bm{\xi}_k(t_1,t))  \varphi_{k}[\mathfrak{u}](t_1,t) \chi(| t_1- t_1^{\bm{x}} |/r) dt_1.
\end{equation}
The evaluation of $\mathcal{J}_{k}(t;\bm{x})$ poses difficulties owing to singularity present in the
integrand at $t=t_{2}^{\bm{x}}$ and the {\em near singularity} in the vicinity of $t_{2}^{\bm{x}}$. In order to circumvent these difficulties, we use
a change of variable $t_1 = t_{1}^{\bm{x}}+\varrho(\tau)$, where the smooth invertible odd function $\varrho(\tau)$ satisfies
\begin{equation}  \label{POLYCOB}
  \frac{d^{m}\varrho(\tau)}{d\tau^m}  \Big\vert_{\tau=0} = 0 \hspace{3mm} \text{for} \hspace{3mm} m=0,...,M,
\end{equation}
(for example, $\varrho(\tau) = \tau^{M+1}$ for even, non-negative integer $M$).
We now change the integration variable in (\ref{eq:-I}) to obtain
\begin{equation} \label{eq:-I-inTau}
 \mathcal{J}_{k}(t;\bm{x}) = \int \limits_{-\varrho^{-1}(r)}^{\varrho^{-1}(r)}G(\bm{x},\bm{\xi}_k(t_1^{\bm{x}}+\varrho(\tau),t)) \varphi_{k}[\mathfrak{u}](t_1^{\bm{x}}+\varrho(\tau),t) \chi(\varrho(\tau)/r) \varrho'(\tau) \,d\tau.
\end{equation}
%where $\tau_{0} = \varrho^{-1}(r_{1})$. 
This change of variable renders the integrand $M$ times differentiable in $\tau$ that are also uniformly bounded. 
In addition, in view of the factor $\chi(\varrho(\tau)/r)$, the integrand in equation (\ref{eq:-I-inTau}) vanishes on the boundary of the integration interval i.e. $\tau = \pm \varrho^{-1}(r)$, together with all of its derivatives. Thus the integrand in integral (\ref{eq:-I-inTau}) 
is a smooth and periodic function which can be integrated to high-order by means of Trapezoidal
rule. Again, this change in a variable produces a set of quadrature point in $\tau$ that
do not coincide with the computational grid on the integration patch $\mathcal{P}_{k}$ requiring an interpolation strategy for which we utilize the FFT refined polynomial interpolation \cite{bruno2001fast}.

The final step in high-order approximation of $ \mathcal{K}_{k}^{sing}[\mathfrak{u}](\bm{x}) $ 
relates to the the computation of $t_2$-integral in (\ref{eq:-VolInt-Adj-Sing-1}). The main difficulty, here,
is encountered in the form of a jump discontinuity in the $t$-derivative of 
$ \mathcal{J}_{k}(t;\bm{x})$  at $t=t_{2}^{\bm{x}}$ \cite{anand2006efficient}. 
%One can see this as a consequence of jump relations in the single-layer potential.
To work around this, we split the integral in (\ref{eq:-VolInt-Adj-Sing-1}) as
\begin{equation}  \label{eq:-VolInt-Adj-Sing-1-Split}
  \mathcal{K}_{k}^{sing}[\mathfrak{u}](\bm{x}) = \int \limits_{0}^{t_{2}^{\bm{x}}} \mathcal{J}_{k}(t;\bm{x}) dt+\int \limits_{t_{2}^{\bm{x}}}^{1} \mathcal{J}_{k}(t;\bm{x}) dt,
\end{equation}
where both integrands, in principle, can be approximated to high-order by means of $Q$-point Newton-Cotes quadrature. 
This, however, presents a practical difficulty in the form of requiring at least $Q$-equidistant grid points in $[0,t_{2}^{\bm{x}}]$, and $[t_{2}^{\bm{x}},1]$,
that, of course, is not available when $t_{2}^{\bm{x}}$ is close to either $0$ or $1$.
The implementation of a $Q$-point quadrature, in such cases, requires values of $\mathcal{J}_{k}(t;\bm{x})$ at points other than the original grid points.
The direct interpolation of $\mathcal{J}_{k}(t;\bm{x})$, however, is neither accurate,
(on account of non-smoothness of  $\mathcal{J}_{k}(t;\bm{x})$ at $t =t_{2}^{\bm{x}}$), nor efficient (because of its dependence on $t_{2}^{\bm{x}}$).
To avoid this expensive computation, we further split (\ref{eq:-VolInt-Adj-Sing-1-Split}) 
as follows:
\begin{eqnarray}  \label{eq:-VolInt-Adj-Sing-MultiSplit}
   \mathcal{K}_{k}^{sing}[\mathfrak{u}](\bm{x}) &=&€Ž\sum_{l_1=1}^{L_1} \int_{t_{2}^{\bm{x}}-l_1(Q-1)h_{2}}^{t_{2}^{\bm{x}}-(l_1-1)(Q-1)h_{2}}\mathcal{J}_{k}(t;\bm{x}) dt  +
  \sum_{i_1=1}^{I_1}\int_{t_{2}^{\bm{x}}-L_1(Q-1)h_{2} -i_1h_{2}}^{t_{2}^{\bm{x}}-L_1(Q-1)h_{2}-(i_1-1)h_{2}}\mathcal{J}_{k}(t;\bm{x}) dt \nonumber\\
  &+& \sum_{l_2=1}^{L_2}\int_{t_{2}^{\bm{x}}+(l_2+1)(Q-1)h_{2}}^{t_{2}^{\bm{x}}+l_2(Q-1)h_{2}}\mathcal{J}_{k}(t;\bm{x}) dt
  +\sum_{i_2=1}^{I_2}\int_{t_{2}^{\bm{x}}+L_2(Q-1)h_{2}+(i_2-1)h_{2}}^{t_{2}^{\bm{x}}+L_2(Q-1)h_{2}+i_2h_{2}} \mathcal{J}_{k}(t;\bm{x}) dt \nonumber \\
  &+&\int_{0}^{t_{2}^{\bm{x}}-L_1(Q-1)h_{2}-I_1h_{2}}\mathcal{J}_{k}(t;\bm{x}) dt + \int_{t_{2}^{\bm{x}}+L_2(Q-1)h_{2}+I_2h_{2}}^{1} \mathcal{J}_{k}(t;\bm{x}) dt,
\end{eqnarray} 
where
$h_{2} = 1/N_2$ denotes the mesh size in $t_2$ direction and 
numbers $L_{1},L_{2},I_{1},I_2$ are obtained as
\begin{align*}
  L_{1} &= \left[ \frac{t_{2}^{\bm{x}}}{(Q-1)h_{2}}\right],\hspace{2mm} L_{2} = \left[ \frac{1-t_{2}^{\bm{x}}}{(Q-1)h_{t}}\right] ,\nonumber\\
  I_{1} &=\left[ \frac{t_{2}^{\bm{x}}-L_{1}(Q-1)h_{2}}{h_{2}}\right],   
  \hspace{2mm}I_{2} =\left[ \frac{(1-t_{2}^{\bm{x}})-L_{2}(Q-1)h_{2}}{h_{2}}\right],
\end{align*}
with $[r]$ denoting the largest integer less than or equal to the real number $r$.

%However, implementation of $Q$-point composite Newton-Cotes quadrature, 
%say $\mathcal{N}^{Q}$, 

Note that when $t_2^{\bm{x}}$ lies on one of the parallel grid lines, that is, $t_2^{\bm{x}} = jh_2$ for some $j$, then last two 
integrals in (\ref{eq:-VolInt-Adj-Sing-MultiSplit}) vanish. Moreover, the computation of the second and the fourth set of integrals therein,
requires at the most $2\times(Q-1)\times(Q-1)$ additional grid points, first $(Q-1)\times(Q-1)$ in the vicinity of $0$ and the remaining in the vicinity of $1$.
%if  instead of breaking integral in Eq.(\ref{eq:-VolInt-Adj-Sing-1}) in 
%two pieces, break it in multiple pieces as 
%at least $Q$-equidistant grid points in $[0,t_{2}^{\bm{x}}]$, and $[t_{2}^{\bm{x}},1]$, that, of course, is not always available,
%particularly when $t_{2}^{\bm{x}}$ close to either $0$ or $1$
%Clearly, multiple splitting of integral in Eq.(\ref{eq:-VolInt-Adj-Sing-1}) solve our both the problem: on one hand it remove the singularity in
%transverse integration and on the other hand it is computationally efficient. 
%In detail, each integrals in Eq.(\ref{eq:-VolInt-Adj-Sing-MultiSplit}) has 
%smooth integrand and it's evaluation via $Q$-point composite Newton-Cotes quadrature, requires function values
%utmost $2\times(Q-1)\times(Q-1)$ additional grid points, first $(Q-1)\times(Q-1)$ in the vicinity of $0$ and the remaining
%in the vicinity of $1$. 
%Furthermore, for any target point $\bm{x}((u_{0},t_{0}))$, the necessary extra grid point required for the 
%integration in Eq.(\ref{eq:-VolInt-Adj-Sing-MultiSplit}) will be the subset of these precomputed $2\times(Q-1)\times(Q-1)$ additional grid points.
%Thus, 
%%to approximate integrals in Eq.(\ref{eq:-VolInt-Adj-Sing-1}) via $Q$-point composite Newton-Cotes quadrature 
%we need additional
%values 
The values, $\mathcal{J}_{k}(t;\bm{x})$ at these additional grid points are obtained by, first interpolating the smooth 
density $\varphi_{k}[\mathfrak{u}]$ at these extra grid lines followed by the integration in (\ref{eq:-I-inTau}). 
%For density $\varphi_{\mathrm{k}}$ interpolation again we have use Fourier
%refined polynomial interpolation. In conclusion, we obtain high-order quadrature for evaluation of singular integral in Eq.(\ref{eq:-VolInt-Adj-Sing-1}).

\begin{remark}
It is clear that the strategy described above is valid when the target point $\bm{x} \in \mathfrak{T}$ happens to lie on one of 
the parallel grid lines. Obviously, the overall computational scheme does require us to compute $\mathcal{K}_{k}^{sing}[\mathfrak{u}](\bm{x})$ 
for target point $\bm{x}$ that do not coincide with any of the grid lines. In particular, this happens when the integration on a boundary patch 
corresponds to a target point coming from an interior patch. To deal with such 
interactions between boundary and interior patches, we first precompute $\mathcal{K}_{k}[\mathfrak{u}](\bm{x})$ at all boundary discretization points
%$\bm{x} \in \cup_{k\in \mathcal{I}_{B}}{\bm{\xi}}_{k}^{-1}(t_{1}^{k},t_{2}^{k})$ 
$\bm{x} \in \cup_{k\in \mathcal{I}_{B}} \mathfrak{T}_k$
and then set up an FFT-refined polynomial interpolation 
schemes in the boundary region $\cup_{k \in \mathcal{I}_B} \mathcal{P}_k$ making possible accurate evaluations of $\mathcal{K}_{k}[\mathfrak{u}]$ at off-grid lines.
\end{remark}

\begin{figure}[h]
\begin{center}
 \includegraphics[clip=true, trim= 12 6 12 16, scale=0.4]{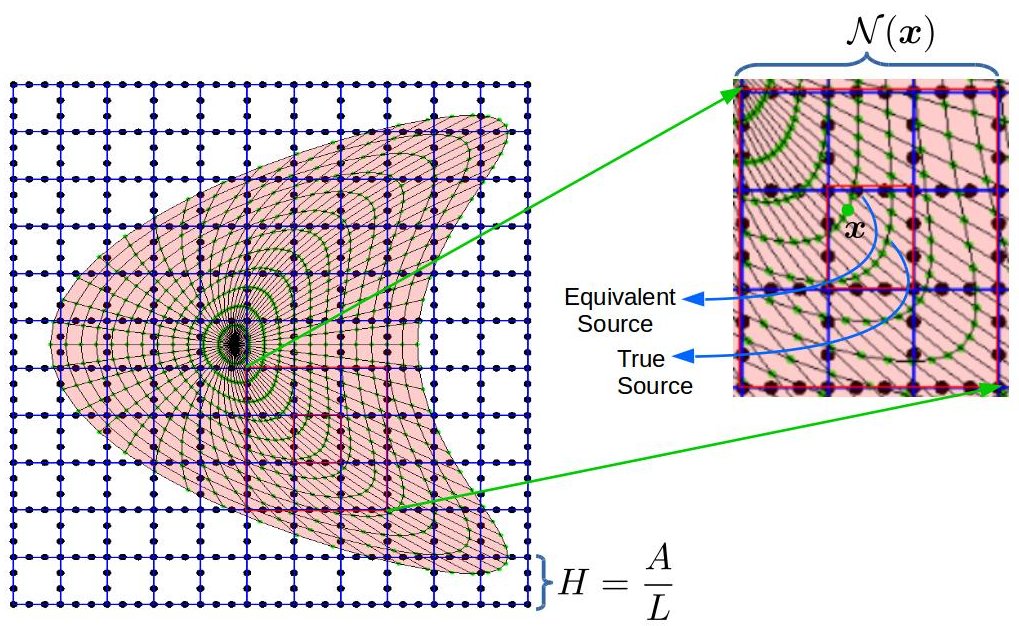}  % l b r t
 \caption{An illustration of the locations $\bm{x}_{ij,\ell}$ of equivalent sources (points displayed in blue color)
	  on the parallel faces of cells $\mathcal{C}_{ij}$. The image also depicts, for a target point $\bm{x}$, the set $\mathcal{N}(\bm{x})$ that defines the adjacency. 
	  %The grid-points within these boxes represent true sources.
	  } \label{fig:-accel}
\end{center}
\end{figure}

 \subsection{Non-singular integration} \label{subsec:-Non-Sing}

%In this section, we describe our acceleration strategy, for the evaluation of
%nonsingular nonadjacent interaction in a fast way. 

As mentioned in the beginning of this section, for the fast computation of non-singular integrals
\begin{align}\label{eq:-VolInt-NonAdj-Reg}
 \mathcal{K}^{reg}[\mathfrak{u}](\bm{x}) &= \left\{ \left( \sum\limits_{k \in \mathcal{I}_B} + \sum\limits_{k \in \mathcal{I}_I} \right)  \int _{\mathcal{P}_{k} \cap\mathcal{N}(\bm{x})^c} \right\}
 G_{\kappa}(\bm{x},\bm{y}) m(\bm{y})\mathfrak{u}(\bm{y})\omega_{k}(\bm{y})\,d\bm{y},
\end{align}
that appear in (\ref{eq:-VolInt-Adj-NonAdj}), we rely on an equivalent sources approximation strategy \cite{bruno2001fast} %, which allows for the use of FFT 
to improve the overall computational efficiency. While the equivalent source methodology is given for the three dimensional Helmholtz kernel in \cite{bruno2001fast}, we can readily adapt those ideas to obtain a variant for the two dimensional case, as discussed below.

%We start by noting that, each one of the volume discretization points $\bm{x} \in  \mathfrak{T}$, that we refer to as a ``true source", lies within one of the 
%cells $\mathcal{C}_{ij}$. % and most of the cells intersect the scatterer $\Omega$.
%As we already pointed out, 
Indeed, we can accurately approximate $\mathcal{K}^{reg}[\mathfrak{u}](\bm{x}_q)$ at each grid point $\bm{x}_q \in \mathfrak{T}$ at an $O(N^2)$ cost by employing a high order quadrature, say, given by
\begin{equation} \label{eq:-disctete-conv}
\mathcal{K}_{a}^{reg}[\mathfrak{u}](\bm{x}_q) = \sum_{\bm{y}_{\ell} \in \mathfrak{T}\setminus \mathcal{N}(\bm{x}_q)} w_{\ell}G_{\kappa}(\bm{x}_q,\bm{y}_{\ell}) m(\bm{y}_{\ell})\mathfrak{u}(\bm{y}_{\ell})\omega_{k}(\bm{y}_{\ell}).
\end{equation}
Noticing that (\ref{eq:-disctete-conv}) is a discrete convolution whose contributing sources $\bm{y}_{\ell}$ are located somewhat irregularly in $\Omega$, 
the acceleration strategy seeks to replace them by a certain set of ``equivalent sources" placed on a Cartesian grid, that produce an accurate approximation
for the convolution, to facilitate the use of FFT for computing the discrete convolution at an $O(N \log N)$ cost.

More precisely, for each cell $\mathcal{C}_{ij}$, if $\mathcal{K}_{a,ij}^{reg}[\mathfrak{u}](\bm{x})$ denotes the quantity
\[
\mathcal{K}_{a,ij}^{reg}[\mathfrak{u}](\bm{x}) = \sum_{\bm{y}_{\ell} \in \mathcal{C}_{ij}} w_{\ell}G_{\kappa}(\bm{x},\bm{y}_{\ell}) m(\bm{y}_{\ell})\mathfrak{u}(\bm{y}_{\ell})\omega_{k}(\bm{y}_{\ell}),
\]
then we seek constants $\sigma^{(m)}_{ij,\ell}$ and $\sigma^{(d)}_{ij,\ell}$ defining the approximating quantity
\begin{equation} \label{equiv-conv-loc}
 \mathcal{K}^{reg,eq}_{a,ij}[\mathfrak{u}](\bm{x}) = \sum \limits_{\mathrm{\ell}=1}^{N^{eq}}\left( \sigma^{(m)}_{ij,\ell}G(\bm{x},\bm{x}_{ij,\ell})+
 \sigma^{(d)}_{ij,\ell}\frac{\partial G(\bm{x},\bm{x}_{ij,\ell})}{\partial \bm{\nu}(\bm{x}_{ij,\ell})}\right).
\end{equation}
such that $$\sum_{q = 1}^{N^{coll}} \left( \mathcal{K}_{a,ij}^{reg}[\mathfrak{u}](\bm{x}_q) - \mathcal{K}_{a,ij}^{reg,eq}[\mathfrak{u}](\bm{x}_q) \right)^2$$ is
minimized for a fixed $N^{coll}$ number of evaluation points $\bm{x}_q$ on the boundary of the set $\mathcal{N}(\bm{x})$ corresponding to the cell $\mathcal{C}_{ij}$. %defined in (\ref{setdef:-Adj})
The locations $\bm{x}_{ij,\ell}$ of equivalent sources in (\ref{equiv-conv-loc}) are equidistant points placed on the two opposite and parallel sides of the cell
% that cell in such a way that the field produced by these equivalent sources accurately approximate the field produced by
%true sources %within the cell $\mathcal{C}_{ij}$ 
%at all point in the space that are at some distance away from the cell $\mathcal{C}_{ij}$. 
%These equivalent sources
%can then be used to compute the inte can be computed efficiently by means of FFTs and then easily transformed back into the scatterer 
%$\Omega$. More precisely, For each cell $c_{ij}$, we seek to replace
% the field produced by the scattering region contained with in cell $c_{ij}$,
%with the field produced by two independent sets of equivalent sources $\mathcal{F}^{l}_{ij},l=1,2$
%, comprising of both acoustical monopoles of intensity $\sigma^{(m),l}_{ij,\mathrm{n}}(\sigma^{(m),l}_{ij,k}G(\bm{x},\bm{x}^{l}_{ij,\mathrm{n}}))$ 
%and acoustical dipoles of intensity $\sigma^{(d),l}_{ij,k}(\sigma^{(d),l}_{ij,k}\frac{\partial G(\bm{x},\bm{x}^{l}_{ij,k})}{\partial y_{l}})$
%placed on an equipoised Cartesian grids $\bm{x}^{l}_{ij,\mathrm{N}},k=1,\cdots,N^{\text{equiv}}$  on the two parallel faces of cell $c_{ij}$
%at all point in the space lying out side to the adjacent set $\mathcal{N}_{x}$, 
as shown in Figure \ref{fig:-accel}.
It is known that one can achieve this approximation up to the prescribed
accuracy $\mathcal{O}(\epsilon)$, provided the number $N^{eq}$ is chosen as \cite{anand2006efficient}
$$ N^{eq} = \max \left \{ \frac{\kappa A}{L},\frac{\log(\epsilon)}{2\log(3/\sqrt{2})}\right \}.$$
%The field generated by these equivalent sources can be obtained as sum of the discrete representation of single and double layer
%potentials of the form
%\begin{equation}
% \mathcal{K}^{l,\text{eq}}_{c_{ij},\text{na}}[\mathfrak{u}](\bm{x}) = \sum \limits_{\mathrm{n}=1}^{\frac{1}{2}N^{\text{equiv}}}\left( \sigma^{(m),l}_{ij,k}G(\bm{x},\bm{x}^{l}_{ij,\mathrm{k}})+
% \sigma^{(d),l}_{ij,\mathrm{n}}\frac{\partial G(\bm{x},\bm{y})}{\partial y_{l}}\Big|_{\bm{y} =\bm{x}^{l}_{ij,\mathrm{n}} }\right).
%\end{equation}
%The unknown intensities $\sigma^{(m),l}_{ij,k}$ and $\sigma^{(d),l}_{ij,k}$ is obtained by minimizing the vector formed by the differences
%$\left(\mathcal{K}^{l,\text{eq}}_{c_{ij},\text{na}}[\mathfrak{u}](\bm{x}) - \mathcal{K}^{l,\text{true}}_{c_{ij},\text{na}}[\mathfrak{u}](\bm{x})\right)$,
%where $\mathcal{K}^{l,\text{true}}_{c_{ij},\text{na}}[\mathfrak{u}](\bm{x})$ stands for field produced by truce sources within cell $c_{ij}$,
% in the mean square norm as $\bm{x}$ varies over a number $n^{\text{coll}}$ collocation points lying on the
% boundary of adjacent set $\mathcal{N}_{x}$. 
The unknown constants $\sigma^{(m)}_{ij,\ell}$ and $\sigma^{(d)}_{ij,\ell}$ are obtained as solution
to the overdetermined linear system ${\bf{A}}  {\bm{\sigma}}^{ij}= \bm{b}$. It is important to note that, as the geometry is identical for each cell $\mathcal{C}_{ij}$,
the $QR$ factorization of the above matrix $\bf{A}$ need only be computed once and saved for repeated use.
For numerical stability of the least-squares solver, we choose $N^{coll}=  4N^{eq}$. A straightforward counting show that this process of 
equivalent source computation for cells requires $O(N^{3/2}/L^3)+O(N^{3/2}/L)$ operations in total.

% as a least square solution of the overdetermined linear system ${\bf{A}}  {\bm{\sigma}}^{ij}= \bm{b}$.
%  
%Having done intensities calculations on uniform Cartesian grids $\bm{x}^{l}_{ij,\mathrm{n}}$, for each $l=1,2$, the nonadjacent field at these point can be obtained as
%  \begin{equation}
%  \mathcal{K}^{l,\text{eq}}_{\text{na}}[\mathfrak{u}](\bm{x}) = \sum \limits_{i,j=1}^{L}\mathcal{K}^{l,\text{eq}}_{c_{ij},\text{na}}[\mathfrak{u}](\bm{x})  
%   -\sum \limits_{\left\{i,j|c_{ij}\subset \mathcal{N}_{x} \right\}}\mathcal{K}^{l,\text{eq}}_{c_{ij},\text{na}}[\mathfrak{u}](\bm{x}).
% \end{equation}
Clearly, for $\bm{x} \in \mathcal{C}_{ij}$, the computation of 
\[
 \mathcal{K}^{reg,eq}_{a}[\mathfrak{u}](\bm{x}) = \sum_{k = 1}^L \sum_{l = 1}^L  \mathcal{K}^{reg,eq}_{a,kl}[\mathfrak{u}](\bm{x}) - \sum_{k = i-1}^{i+1} \sum_{l = j-1}^{j+1}  \mathcal{K}^{reg,eq}_{a,kl}[\mathfrak{u}](\bm{x})
\]
as an approximation to $ \mathcal{K}^{reg}_{a}[\mathfrak{u}](\bm{x})$ at all points on the grid can now be obtained  by means of FFTs with a computational cost of 
$O\left(LN^{1/2}\log (LN^{1/2})\right)$.  
The values obtained in this manner provide 
accurate approximations for non-adjacent non-singular interactions
at points $\bm{x}_{ij,\ell}$ on the boundary of cells $\mathcal{C}_{ij}$.
Finally, to obtain these values
%the approximation of nonadjacent nonsingular interactions, i.e.,
%\begin{align}
%  \mathcal{K}^{\text{true}}_{\text{na}}[\mathfrak{u}](\bm{x}) &=\underset{\bm{y}\in \mathcal{P}_{\mathrm{k}} \cap\mathcal{N}_{x}^{c}}{\sum\limits_{\mathrm{k} \in \mathcal{I}_{B}}}\int \limits_{0}^{1} \int \limits_{0}^{1}
%  G(\bm{x},\bm{y}(u,t))  \varphi_{\mathrm{k}}(u,t) du dt  \nonumber  \\
%   &+\underset{\bm{y}\in \mathcal{P}_{\mathrm{k}} \cap\mathcal{N}_{x}^{c}}{\sum\limits_{\mathrm{k} \in \mathcal{I}_{I}}}\int \limits_{0}^{1} \int \limits_{0}^{1}
%  G(\bm{x},\bm{y}(u,t))  \varphi_{\mathrm{k}}(u,t) du dt
%\end{align}
at a true source location, say for an $\bm{x} \in \mathcal{C}_{ij} $, we solve the
free space Helmholtz equation within $\mathcal{C}_{ij}$, with Dirichlet boundary boundary data coming from $\mathcal{K}^{req,eq}_{a}[\mathfrak{u}](\bm{x}_{ij,\ell})$.
%The absence of internal resonance ensures that 
To efficiently obtain solutions to these well-posed Dirichlet boundary value problems, %are well-posed where solutions can be obtained by means 
we utilize a discretized plane wave expansion of the form \cite{felsen1994radiation} 
\begin{equation} \label{eq:-PWE}
  \mathcal{K}^{req,eq}_{a}[\mathfrak{u}](\bm{x}) \approx \sum \limits_{\ell=1}^{N^{w}}\gamma_{j}\text{exp}(i\kappa \bm{d}_{\ell}.\bm{x}),
\end{equation}
where the unit vectors $\bm{d}_{j}$  sample the surface of unit disc with sufficient degree of uniformity. 
Use of this approach is motivated by the spectral convergence of the above expansion with respect to the number of unit vectors  $\bm{d}_{\ell}$ used and the fact that 
the wave expansion coefficients $\gamma_{\ell}$  can be obtained as
%we equate expansion (\ref{eq:-PWE}) with $ \mathcal{K}^{l,\text{eq}}_{\text{na}}[\mathfrak{u}](\bm{x})$ 
%in the least square sense at all nonadjacent point located at the the boundary of cell $c_{ij}$. This procedure 
solutions to overdetermined linear systems of the form $ \bm{B \gamma} = \bm{\beta_{ij}}$ where the matrix $ \bm{B}$  remains unchanged for each cell, again, owing to the identical geometry of cells $\mathcal{C}_{ij}$s.
%thus en a single $QR$-decomposition of $\bm{B}$ to be used for multiple solves. 
Thus, $\mathcal{K}^{reg,eq}_{a}[\mathfrak{u}](\bm{x})$ can be evaluated at all true source locations
$\bm{x}\in \mathcal{C}_{ij}$ at a computational cost of $O(N^{3/2}/L^3)+O(N^{3/2}/L)$ and the overall cost 
of evaluation of nonadjacent non-singular interactions, therefore, stands at $O\left(LN^{1/2}\log (LN^{1/2})\right)+O(N^{3/2}/L^3)+O(N^{3/2}/L)$.

\begin{table}[t]
 \begin{center}
\begin{tabular}{c|c|c|c|c|c|c}\hline 
Grid Size &  Unknown & Iteration & \multicolumn{2}{c|}{$L^2$} & \multicolumn{2}{c}{ $L^{\infty}$ } \\
 \cline{4-7}
 &   & Number & $\varepsilon_{2}$ & Order  & $\varepsilon_{\infty}$ &  Order \\ 
\hline
$2\times3\times9+ 1\times9\times$9 & 135 &  10 & 7.04e-01 & - & 6.21e-01 & -  \\
\hline
$2\times5\times17+ 1\times17\times$17 & 459 &  11&2.68e-01 & 1.39e+00&2.73e-01 & 1.18e+00 \\
\hline
$2\times9\times33+ 1\times33\times$33 & 1683 &  12 & 3.92e-02 & 2.78e+00&4.77e-02 & 2.52e+00 \\
\hline
$2\times17\times65+ 1\times65\times$65 & 6435 & 15 &3.38e-03 & 3.54e+00&4.14e-03 & 3.52e+00 \\
\hline
$2\times33\times129+ 1\times129\times$129 & 25155 &  18&1.47e-04 & 4.52e+00&2.37e-04 & 4.13e+00 \\
\hline
\end{tabular}
\caption{Convergence study: plane wave scattering by a disc of acoustic size $\kappa a =10$ and refractive index $n = \sqrt{2}$ when $3$-point Newton-Cotes quadrature employed for $t_2$-integration over boundary patches.}  \label{table:-Disc-Q3-K5}
 \end{center} 
 \begin{center}
\begin{tabular}{c|c|c|c|c|c|c} \hline
Grid Size &  Unknown & Iteration &\multicolumn{2}{c|}{ $L^2$ } & \multicolumn{2}{c}{ $L^{\infty}$ } \\ 
\cline{4-7}
 &   & Number & $\varepsilon_{2}$ & Order  & $\varepsilon_{\infty}$ &  Order \\ 
\hline
$2\times5\times9+ 1\times9\times9$ & 171 &  10 & 7.09e-01 &- &5.96e-01& -  \\
\hline
$2\times9\times17+ 1\times17\times17$ & 595 &  11&1.22e-01&  2.54e+00 & 1.27e-01 & 2.24e+00 \\
\hline
$2\times17\times33+ 1\times33\times33$ & 2211 &  15& 8.78e-03 & 3.80e+00 & 1.03e-02 & 3.62e+00  \\
\hline
$2\times33\times65+ 1\times65\times65$ & 8515 &  18&2.84e-04& 4.95e+00 & 3.35e-04 & 4.94e+00 \\
\hline
$2\times65\times129+ 1\times129\times129$ & 33411 &  23&7.76e-06 & 5.19e+00  & 8.53e-06 & 5.30e+00 \\
\hline
\end{tabular}
\caption{Convergence study: plane wave scattering by a disc of acoustic size $\kappa a =10$ and refractive index $n = \sqrt{2}$ when $5$-point Newton-Cotes quadrature employed for $t_2$-integration over boundary patches.}  \label{table:-Disc-Q5-K5}
\end{center}

 \end{table}

\subsection{Numerical Results}

In this section, we present a series of numerical examples to exemplify the effectiveness of the
two dimensional algorithm discussed in previous subsections in terms of its numerical accuracy and computational efficiency. Be begin by noting that, in all convergence tables given below, we use
$p_{1}\times n_{1}\times n_{2}+q_1\times m_1\times m_2$ to announce the underlying computational grid implying that a $p_{1}$ number of boundary patches, each with $n_{1}\times n_{2}$ discretization points, 
and a $q_{2}$ number of interior patches, each with $m_{1}\times m_{2}$ points, constituted the Nystr\"{o}m grid $\mathfrak{T}$. 
The relative error (in the near field) 
%$\varepsilon_{\infty}$, $\varepsilon_{2}$ 
reported in tables are obtained as 
\begin{align*}
 \varepsilon_{\infty} &= \frac{ \underset{1\leq i \leq N}{\max} \left|\mathfrak{u}^{\text{exact}}(\bm{x}_{i})-\mathfrak{u}^{\text{approx}}(\bm{x}_{i})\right|}{ \underset{1\leq i \leq N}{\max} \left|\mathfrak{u}^{\text{exact}}(\bm{x}_{i}) \right|}, \\
  \varepsilon_{2} &= \left ( \frac{\sum \limits_{i=1}^{N}\left|\mathfrak{u}^{\text{exact}}(\bm{x}_{i})-\mathfrak{u}^{\text{approx}}(\bm{x}_{i})\right|^2}{\sum \limits_{i=1}^{N}\left|\mathfrak{u}^{\text{exact}}(\bm{x}_{i}) \right|^2} \right)^\frac{1}{2}.
\end{align*}
We write ``Order"  to denote the numerical order of convergence of the approximation. %calculated as logarithmic of ratio of the relative error $\varepsilon_{\infty}$, $\varepsilon_{2}$ at successive level of discretization.  
All convergence studies correspond to the scattering of incident plane waves propagating along the $x$-axis. As our approximation scheme for the integral operator is based on the use
of spectrally accurate Trapezoidal rule and $Q$-point Newton-Cotes quadrature, the expected convergence rates are essentially that of Newton-Cotes, that is, either $Q$ or $Q+1$ depending on whether $Q$ is even or odd, provided other parameters, such as the degree of the polynomial in the interpolators, the number $M$ for the change of variable $\varrho$, etc. are chosen favorably.
%Therefore the expected order of convergence of the algorithm is precisely the order of convergence ofcomposite Newton-Cotes quadrature used for transverse integration over boundary patches. 
%Moreover, one can certainly increase
%the order of convergence of the algorithm by using more accurate composite Newton-Cotes quadrature. 
We demonstrate the increase in the convergence rate with respect to increasing $Q$ through two different sets of experiments, one using a $3$-points while other using $5$-points Newton-Cotes quadrature.

\begin{figure}[h!] \label{fig:-Disc-K10}
\begin{center}
\subfigure[Disc]{\includegraphics[clip=true, trim=400 45  350 30, scale=0.27 ]{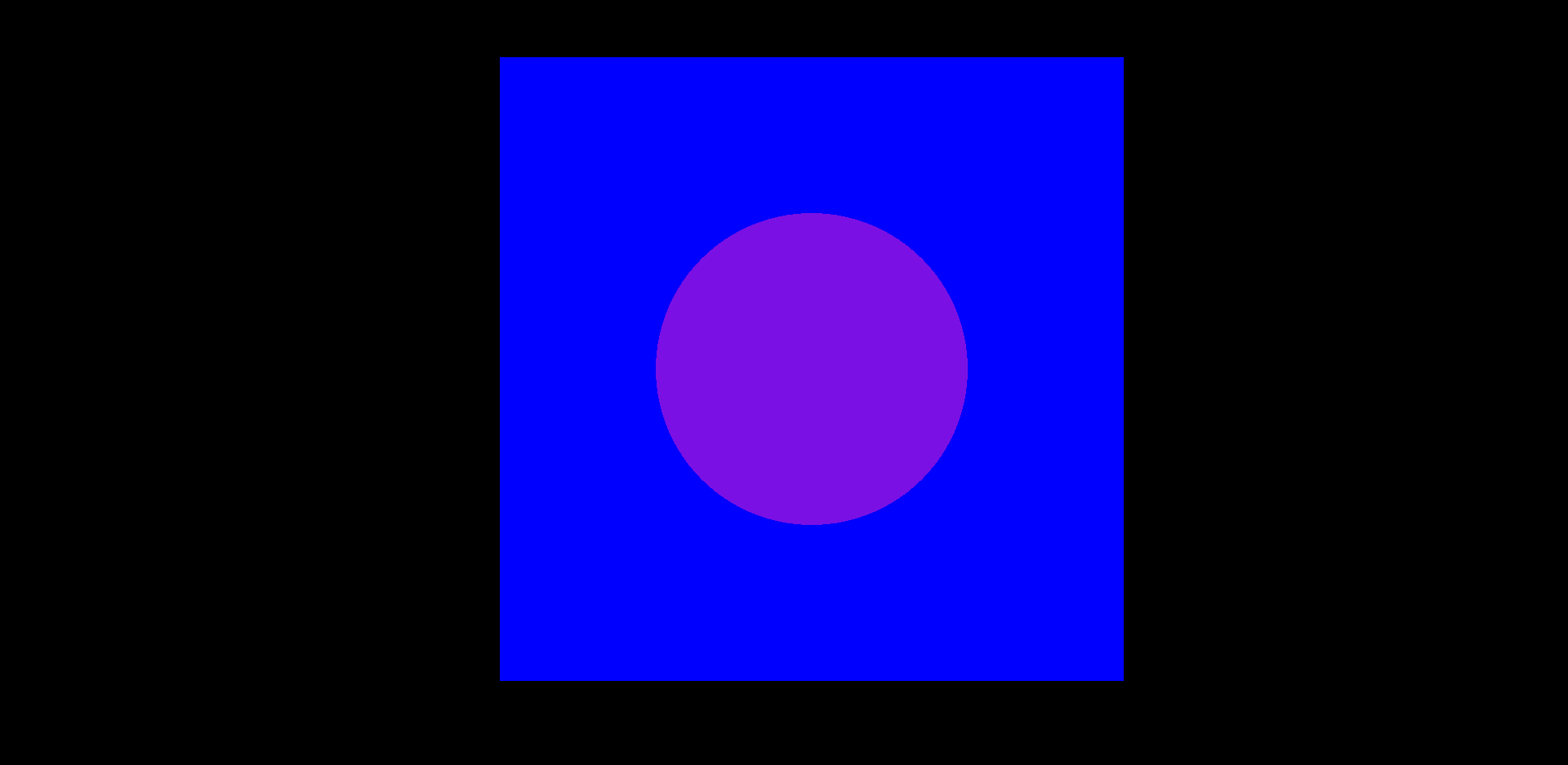}} % l b r t
\hfill
\subfigure[Plane wave incidence, $\Re(\mathfrak{u}^i)$]{\includegraphics[clip=true, trim=400 45  350 30, scale=0.27]{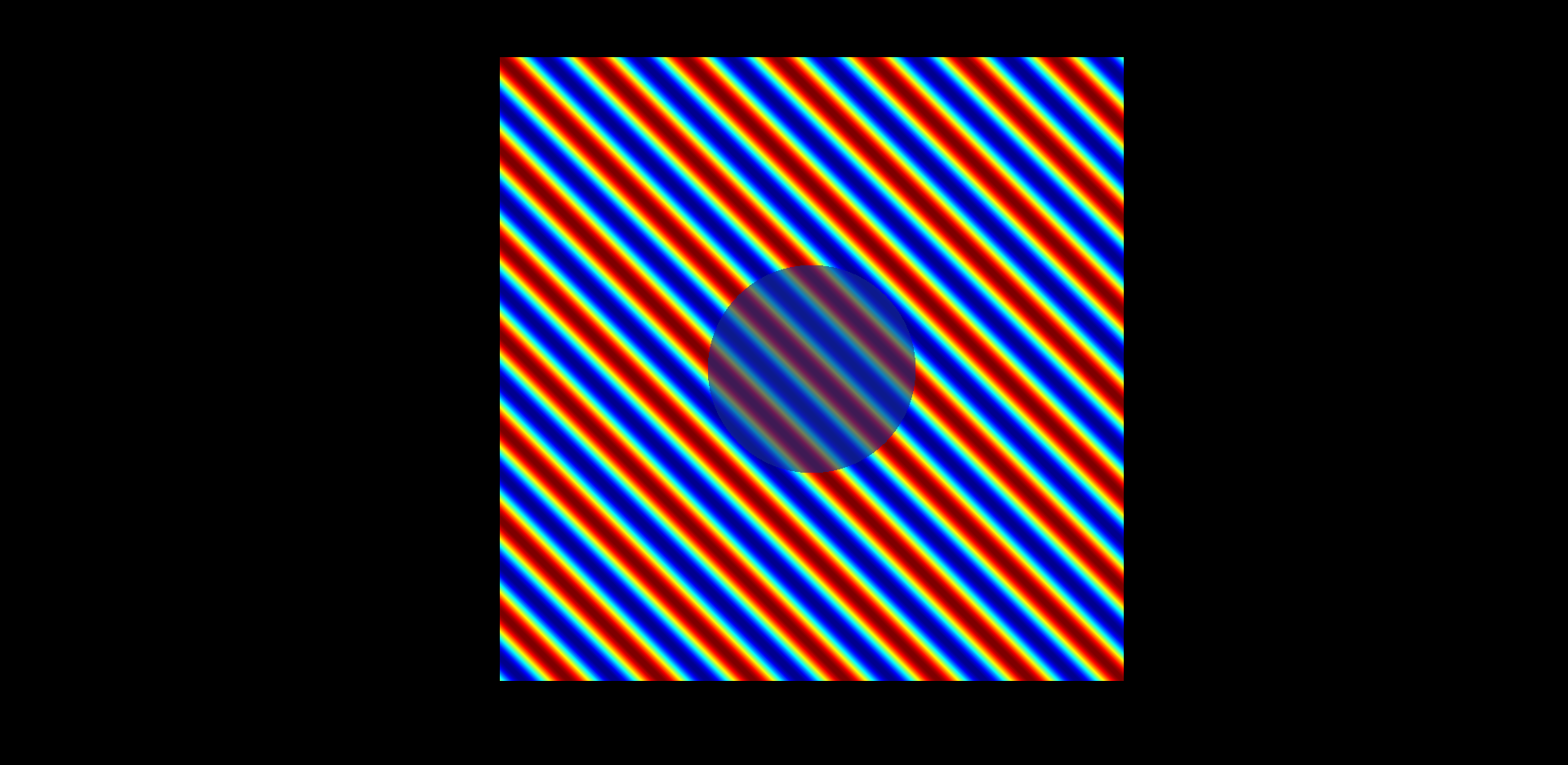}}
\hfill
\subfigure[Real Part of the scattered field, $\Re(\mathfrak{u}^s)$]{\includegraphics[clip=true, trim=400 45  350 30, scale=0.27]{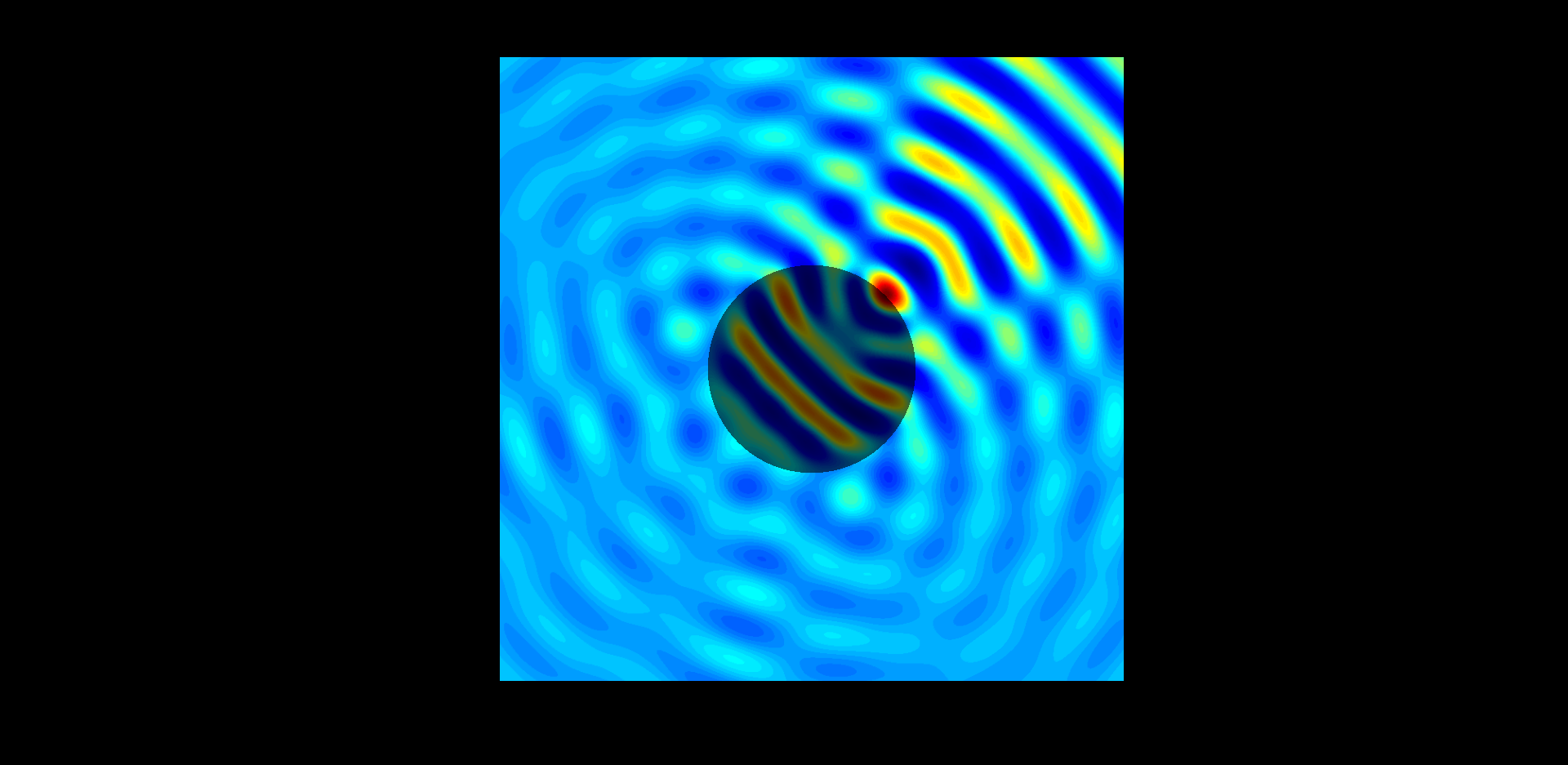}}
\hfill
\subfigure[Real part of the total field, $\Re(\mathfrak{u})$]{\includegraphics[clip=true, trim=220 30  200 10, scale=0.49]{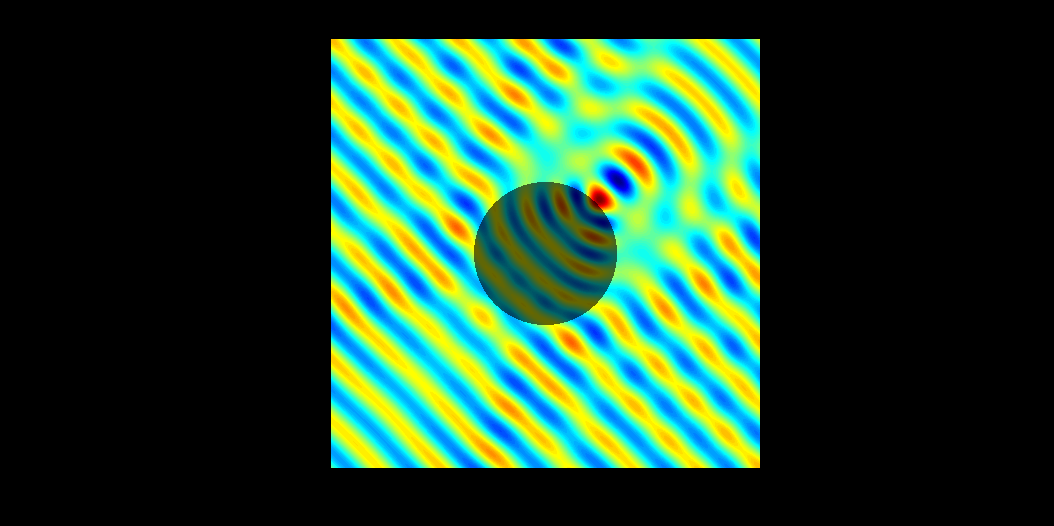}}
\hfill
\caption{Scattering of a plane wave $ \exp (i \kappa \bm{d} \cdot \bm{r})$ with  $\kappa = 10, d = (1/\sqrt{2},1/\sqrt{2})$ by a penetrable unit disc. A computational grid of size $2\times{65} \times{129}+ 1 \times 129 \times 129 $ and $5$-point Newton-Cotes quadrature
	 is used to obtain an accuracy of 0.006\%.}\label{fig:-Disc-K10}
\end{center}
\end{figure}

\begin{exmp}(\textit{Convergence study for a simple scatterer})

 As a first exercise, we compute the scattering by a penetrable disc
 for which a comparison with exact solution is possible.
 For this scattering configuration, we present a convergence study, in Table \ref{table:-Disc-Q3-K5}, when the acoustic size ($\kappa a$, $a$ being the diameter of the inhomogeneity) of the scatterer is $10$ and its refractive index is $n = \sqrt{2}$.
% have conducted three numerical experiment for
%$(1) \kappa a = 4, n =\sqrt{3}$ $(2) \kappa a = 10, n= \sqrt{2}$ $(3) \kappa a = 20, n =\sqrt{2}$,
%where $a$ stands for the diameter of the disc. For $\kappa a = 4$ and $\kappa a = 10 $
The numerical solutions in this study are computed for several discretization levels while employing a
$3$-point Newton-Cotes quadrature for $t_2$-integration in boundary patches. In this case, we clearly see the methodology
achieve the expected rate of convergence. To show an enhanced order of convergence, we repeat the same experiment but use a $5$-point Newton-Cotes quadrature in the computation of the integral operator. The corresponding results are presented in Table \ref{table:-Disc-Q5-K5} where, as expected, we see an increased rate of convergence.

We display, in Figure \ref{fig:-Disc-K10}, the numerical solution obtained using our methodology  
corresponding to an experiment with $ \kappa a = 20$ and $n =\sqrt{2}$
where the incidence wave impinges on the obstacle at $45$ degree angle with positive $x$-axis. This computation required $57$ GMRES iterations to produce an accuracy of 0.006\% in the max norm.
%reduce the residual by a factor of $10^{-5}$.

\end{exmp}

\begin{table}[h!]
\begin{center}
\begin{tabular}{c|c|c|c|c|c} \hline
$\kappa a$ & Grid Size & $\varepsilon_{2}$  & $\varepsilon_{\infty}$  & \multicolumn{2}{c} {Time(s)/Iteration} \\
\cline{5-6}
& & & & accelerated & Non-Accelerated \\
\hline
4 & $2\times 16\times 32 + 1\times 32 \times 32$ & 1.18e-03 &  3.15e-03  & 2.57  & 3.93  \\
\hline
8 & $2\times 32 \times 64 + 1\times 64 \times 64$ & 3.90e-04 & 5.67e-04  & 8.53  & 39.56  \\
\hline
16 & $2\times 64 \times 128 + 1\times 128 \times 128$ & 7.90e-04 & 7.85e-04 & 32.8  & 626.47  \\
\hline
32 & $2\times 128 \times 256 + 1\times 256 \times 256$ & 8.64e-05  & 9.05e-05 & 122.31 & 10785.9\\
\hline
\end{tabular}
\caption{Computational cost for accelerated and non-accelerated algorithm.} \label{table:-effi}
\end{center}
\end{table}

\begin{exmp} (\text{Computational Efficiency})
We next demonstrate a tempered growth in computational cost of the proposed accelerated algorithm while
maintaining a fixed computational error by comparing our approximate integral operator against the continuous operator (\ref{eq:-VolInt}).
Toward this, with every doubling of the computational grid size, 
we double the wave number $\kappa$ to keep the number of points per wavelength unchanged, thus, fixing the accuracy level.
For this set of experiments, we again use a disc scatterer with a constant refractive index $n=2$. The results in Table \ref{table:-effi} 
show that, while for small values of $N$, the cost of non-accelerated algorithm is comparable to that of the accelerated version, 
as problem size increase, there is a substantial gain in terms of computational cost.
% of accelerated methodology  substantially reduced compare to the non-accelerated counterpart. 
For instance, for the scatterer of size $\kappa a = 32$, accelerated computations are $88$ times faster than its non-accelerated counterpart. In particular, the time for accelerated computation in Table \ref{table:-effi} exhibit the growth according to the computational complexity of $O(N \log N)$.

\end{exmp}

\begin{table}[t]
\begin{center}
\begin{tabular}{c c c|c|c c c c c c c c c  |c c c c c c c c c }\hline
& $\kappa H$ & & Equivalent Source/cell & & & &  $\varepsilon_{2}$ & & & &  & & & & $\varepsilon_{\infty}$ & & & & \\
\hline
 & 2 & & $4\times 4 $& & & &  5.10e-05 & & & &  & & & & 1.30e-04 & & & &\\
\hline
 & 4 & &$8\times 4 $& & & & 3.05e-09 & & & &  & & & & 1.94e-08 & & & & \\
\hline
 & 8 & &$13\times 4 $& & & &  3.00e-12 & & & &    & & & & 6.70e-12 & & & &\\
\hline
 & 16 &  & $30\times$ 4 & & & &  1.74e-12 & & & &  & & & &  1.14e-12& & & &\\
   \hline
\end{tabular}
\caption{Accuracy of acceleration when wavenumber $\kappa$ increases but number of point per wavelength is fixed.} \label{table:-accel2}
\end{center}
\end{table}

\begin{exmp} (\textit{Accuracy of non-singular non-adjacent interactions})

\begin{table}[t]
 \begin{center}
\begin{tabular}{c| c c c c c c c c c c | c c c c c c c c c c c } \hline
Equivalent Source/cell & & & & & $\varepsilon_{2}$ & & & & & & & & &  &  $\varepsilon_{\infty}$ & & & & &   \\
\hline
$4 \times 4$ & & &  & & 1.60e-02 & & & & &  & & & & &  4.66e-02 & & & & &   \\
\hline
$6 \times 4$& & & & &  4.14e-05 & & & & &  & & & & &  1.68e-04 & & & & &    \\
\hline
$8 \times 4$& & & & &3.53e-08 & & & & &  & & & &  &  1.85e-07  & & & & &   \\
\hline
$10 \times 4$ & & & & & 1.60e-10& & & & & & & & & &  1.48e-09 & & & & &  \\
\hline
$12 \times 4 $& & & & & 2.86e-12 & & & & & & &  & & &  2.82e-11& & & & &   \\
  \hline
\end{tabular}
\caption{Accuracy of acceleration for fixed wavenumber $\kappa = 8\pi$}\label{table:-accel1}
\end{center}
\end{table}

%As seen in previous examples, the acceleration strategy brings considerable gain in terms of overall computational cost of the methodology. We now 
% investigate
% computed non-singular non-adjacent interactions at true source point
%by means of equivalent sources placed on parallel faces of cell $c_{ij}$. This example 

The previous examples demonstrate the high-order convergence and efficient computational times for the proposed methodology.
We, in particular, noted significant gains for the accelerated scheme in terms of computational cost. On the other hand,
though Table \ref{table:-effi} does show that the accuracy levels of the numerical method are maintained even when frequency of wave oscillations 
increase, provided we ensure a fixed number of discretization points per wavelength, 
an interesting picture emerges when we look approximations coming from the accelerator alone. 
To see this, we compare the values of $\mathcal{K}_{a,ij}^{reg, eq}[\mathfrak{u}]$ against that of $\mathcal{K}_{a,ij}^{reg}[\mathfrak{u}]$ at increasing levels of
frequency while the number of equivalent sources remained unchanged. As seen in Table \ref{table:-accel2}, approximations converge rapidly with increasing wavenumber.
This observation, in turn, is utilized in a more effective load balancing between adjacent and non-adjacent calculations to achieve favorable computational cost.

We also include, in Table \ref{table:-accel1}, a similar study, this time keeping the wavenumber constant ($\kappa = 8 \pi$) while increasing the number of equivalent sources, $N^{eq}$, 
where we again see rapid improvement in accuracy levels. Tables \ref{table:-accel2} and  \ref{table:-accel1} clearly demonstrate
the spectral accuracy of acceleration strategy for small as well as large cell sizes $H$. Indeed, error introduced as a result of acceleration procedure are typically much smaller in comparison with other
sources of error within the overall algorithm, and have no impact on convergence rates, as confirmed by the numerical studies in this section.

\end{exmp}

 \begin{figure}[h!] \label{fig:-Bean-K10}
\begin{center}
\subfigure[Bean]{\includegraphics[clip=true, trim=400 45  350 30, scale=0.28 ]{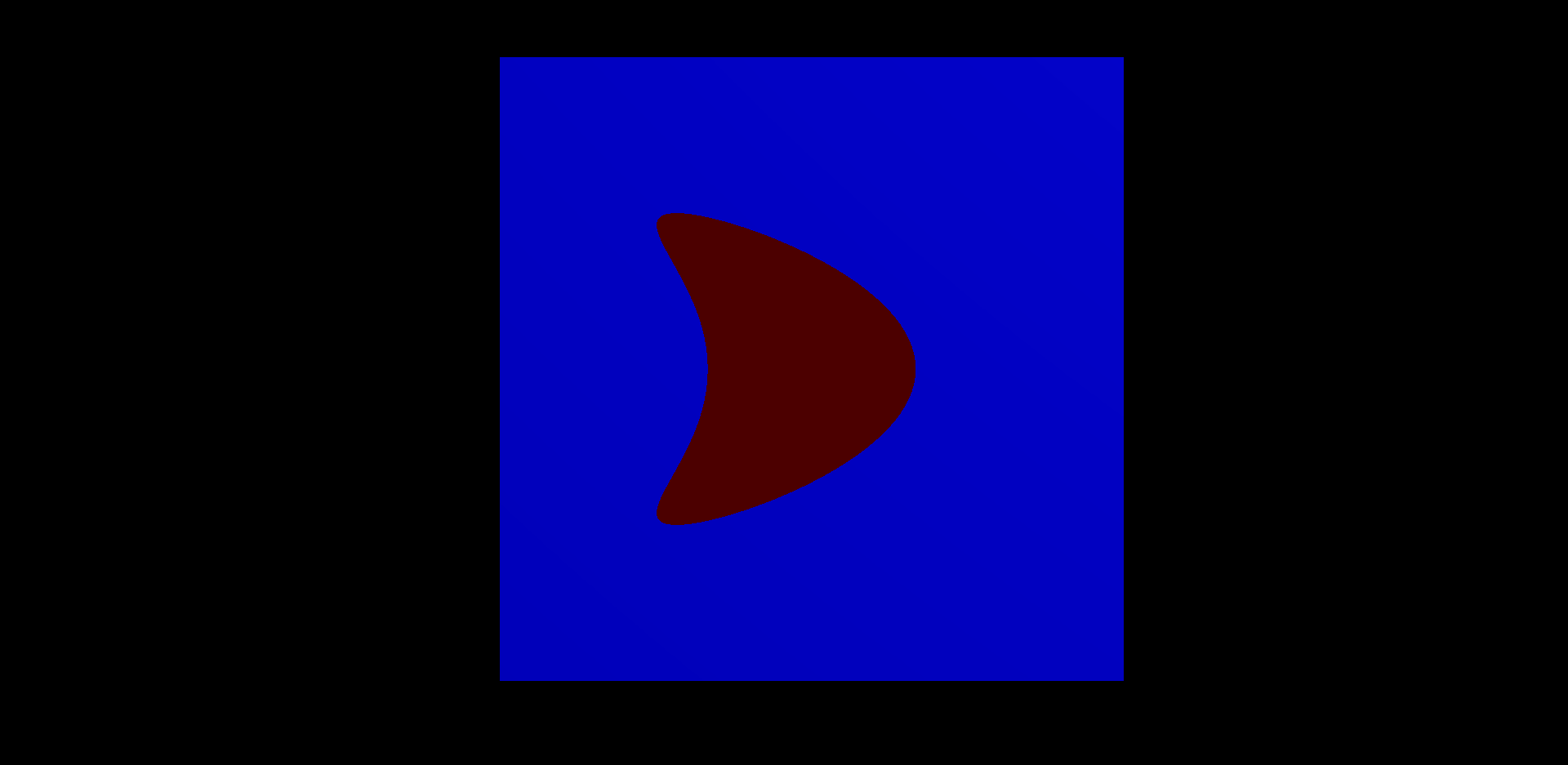}} % l b r t
\hfill
\subfigure[The incident plane wave, $\Re(\mathfrak{u}^i)$]{\includegraphics[clip=true, trim=400 45  350 30, scale=0.28]{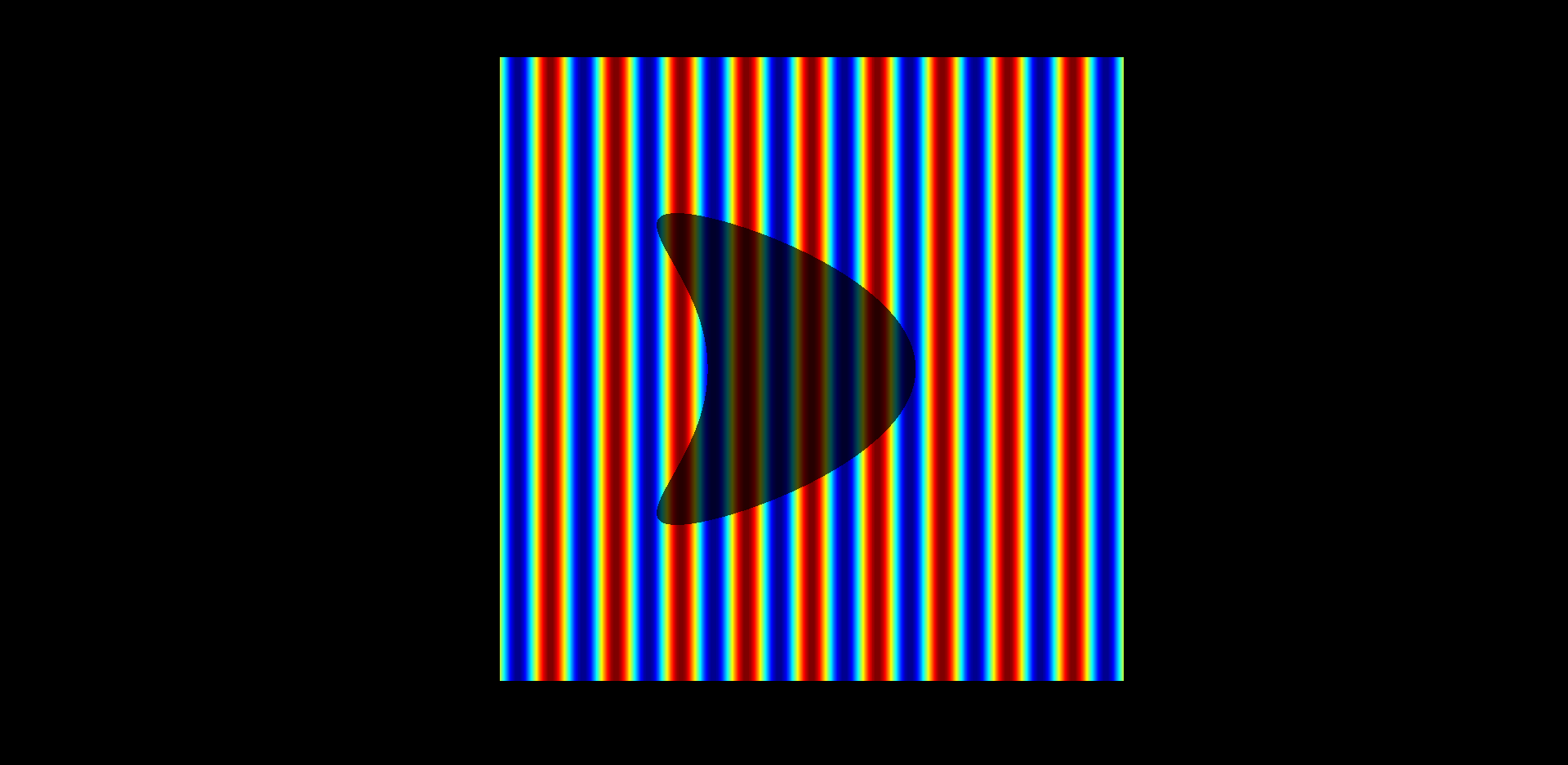}}
\hfill
\subfigure[Real Part of the scattered field, $\Re(\mathfrak{u}^s)$]{\includegraphics[clip=true, trim=400 45  350 30, scale=0.28]{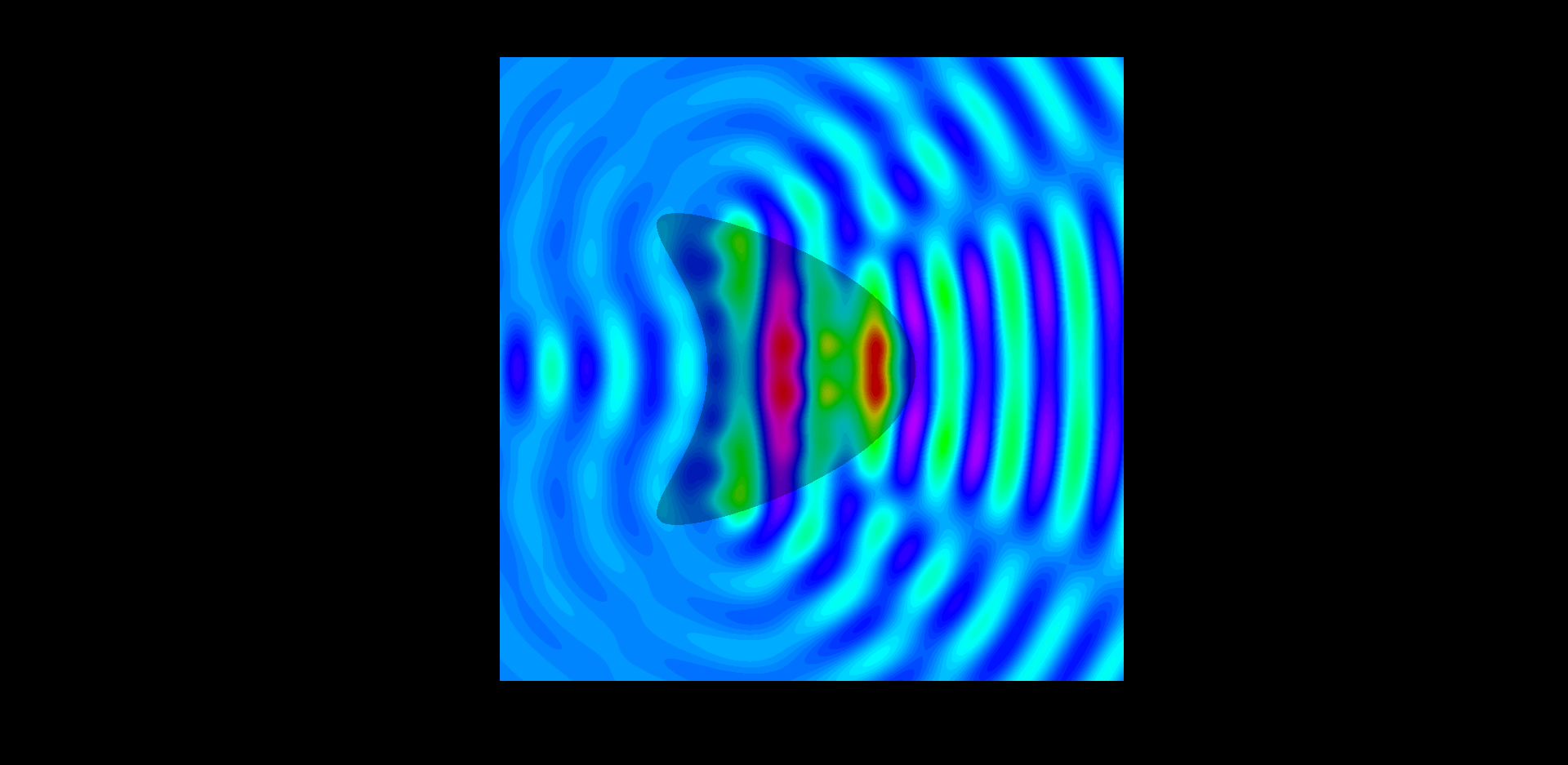}}
\hfill
\subfigure[The real part of the total field, $\Re(\mathfrak{u})$]{\includegraphics[clip=true, trim=400 45  350 30, scale=0.28]{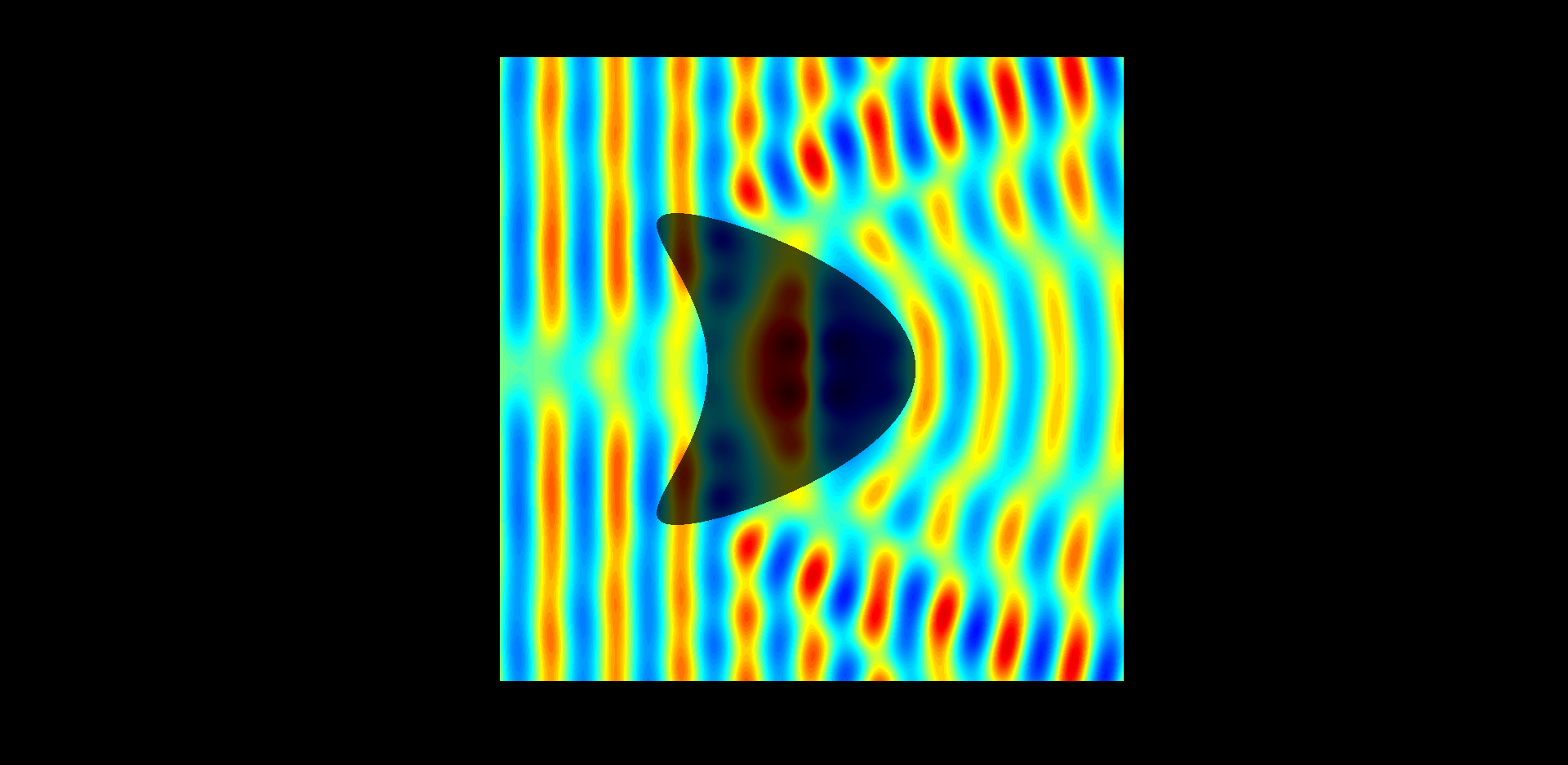}}
\hfill
\caption{Scattering by penetrable bean shape scatterer, At frequency $k= 10, d = (1,0)$, 
	 Grid Size $2\times{65} \times{129}+ 1 \times 129 \times 129 $, Newton-Cotes five point quadrature
	 is used for transverse integration over boundary patches.}\label{fig:-Bean-K10}
\end{center}
\end{figure} 

\begin{exmp}(\textit{Convergence study: complex scattering media})

The methodology presented in this text can, of course, be applied to acoustic scattering calculations from penetrable media with complicated geometries as well as variable material properties.

We begin with an example that demonstrates adaptability and applicability of our algorithm in dealing with scatterers that have relatively complex geometrical description. Toward this, we consider scattering by a penetrable bean shaped scatterer, as shown in Figure \ref{fig:-Bean-K10}. For these scattering computations, we use an incident plane wave $\exp(i \kappa x)$ with $\kappa a = 10$. The refractive index of the medium is taken to be $n = \sqrt{2}$.
We use a numerical solution at a fine grid for convergence study comparisons. 
%therefore in order to obtain relative error computed solution compared against the solution obtained
%with a finer discretization. For this scattering configuration, 
The first sets of computations, for which the results are presented in Table \ref{table:-Bean-K5-Q3},  employ $3$-point Newton-Cotes quadrature
for $t_2$-integration over boundary patches. As expected, we see a convergence rate of order $4$ in the solution.
We repeat the experiment with
$5$-points Newton-Cotes quadrature and report the results in Table \ref{table:-Bean-K5-Q5}, where we again see an enhanced order of convergence. 
 
\begin{table}[t!]
 \begin{center}
\begin{tabular}{c|c|c|c|c|c|c} \hline\label{2}
Grid Size &  Unknown & Iteration &\multicolumn{2}{c|}{ $L^2$ } & \multicolumn{2}{c}{ $L^{\infty}$ } \\
\cline{4-7}
&   & Number & $\varepsilon_{2}$ & Order  & $\varepsilon_{\infty}$ &  Order \\
\hline
$2\times3\times9+ 1\times9\times$9 & 135 & 4&  1.95e+00 & - &1.82e+00 & -\\
\hline
$2\times5\times17+ 1\times17\times$17 & 459 & 11&  3.89e-01 & 2.33e+00&3.84e-01 & 2.25e+00 \\
\hline
$2\times9\times33+ 1\times33\times$33 & 1683 &15 &  3.88e-02 & 3.32e+00&4.07e-02 & 3.24e+00\\
\hline
$2\times17\times65+ 1\times65\times$65 & 6435 &20 &  4.10e-03 & 3.24e+00&4.83e-03 & 3.08e+00\\
\hline
$2\times33\times129+ 1\times129\times$129 & 25155 & 25&  2.28e-04 & 4.17e+00&4.55e-04 & 3.41e+00 \\
\hline
\end{tabular}
\caption{Convergence for the bean shape scatterer with $\kappa =5$ and $n = \sqrt{2}$,
	 when $3$-point Newton-Cotes quadrature
	 employed for boundary patch integration in transverse integration.}\label{table:-Bean-K5-Q3}
\end{center}
\end{table} 

 \begin{table}[b!]
 \begin{center}
\begin{tabular}{c|c|c|c|c|c|c}\hline\label{2}
Grid Size &  Unknown & Iter & \multicolumn{2}{c|}{$L^2$} & \multicolumn{2}{c}{ $L^{\infty}$ } \\
 \cline{4-7}
 &   &  & $\varepsilon_{2}$ & Order  & $\varepsilon_{\infty}$ &  Order \\ 
\hline
$2\times5\times9+ 1\times9\times9$ & 171 & 12 & 1.85e+00& & 1.31e+00 & - \\
\hline
$2\times9\times17+ 1\times17\times17$ & 595 &  16 &1.56e-01& 3.57e+00 & 3.35e-01 & 1.96e+00\\
\hline
$2\times17\times33+ 1\times33\times33$ & 2211 &  18 &1.90e-02 &3.04e+00 & 2.94e-02 & 3.51e+00 \\
\hline
$2\times33\times65+ 1\times65\times65$ & 8515 & 24 & 6.13e-04& 4.95e+00& 6.36e-04 & 5.53e+00 \\
\hline
$2\times65\times129+ 1\times129\times129$ & 33411 & 30  &1.17e-05 & 5.71e+00 & 1.30e-05 & 5.61e+00\\
\hline
\end{tabular}
\caption{Convergence for the bean shape scatterer with $\kappa =5$ and $n = \sqrt{2}$,
	 when $5$-point Newton-Cotes quadrature
	 employed for boundary patch integration in transverse integration.}\label{table:-Bean-K5-Q5}
\end{center}
\end{table} 
 
Of course, our algorithm is not restricted to consideration of constant or piecewise constant material properties. In order to demonstrate this fact, in our final experiment, we carried out scattering computations for 
a variable refractive index $n$ given by
\[
n(\bm{x}) = \sin(\pi x_1)\cos(\pi x_2) \hspace{2mm}\text{ for } \bm{x} = (x_1,x_2) \in \Omega.
\]
A plane wave incidence with $ \kappa =2$ is used and results corresponding to three and five point Newton-Cotes quadrature
 are presented in table (\ref{table:-Bean-K2-Q3}) and (\ref{table:-Bean-K2-Q5}) respectively. Convergence study of these tables
 clearly support that the high-order nature of our algorithm, as expected, continues to hold for scattering media with variable material properties.

\begin{table}[t!]
 \begin{center}
\begin{tabular}{c|c|c|c|c|c|c} \hline
Grid Size &  Unknown & Iteration &\multicolumn{2}{c|}{ $L^2$ } & \multicolumn{2}{c}{ $L^{\infty}$ } \\
\cline{4-7}
&   & Number & $\varepsilon_{2}$ & Order  & $\varepsilon_{\infty}$ &  Order \\
\hline
$2\times3\times9+ 1\times9\times$9 & 135 & 3 &1.01e+00 & - &1.25e+00 & -  \\
\hline
$2\times5\times17+ 1\times17\times$17 & 459 &  4&7.60e-02 & 3.74e+00&9.11e-02 & 3.77e+00 \\
\hline
$2\times9\times33+ 1\times33\times$33 & 1683 & 5 &7.71e-03 & 3.30e+00&8.36e-03 & 3.45e+00 \\
\hline
$2\times17\times65+ 1\times65\times$65 & 6435 & 6 &6.56e-04 & 3.55e+00&6.36e-04 & 3.72e+00\\
\hline
$2\times33\times129+ 1\times129\times$129 & 25155 & 7 &2.02e-05 & 5.02e+00&5.79e-05 & 3.46e+00 \\
\hline
\end{tabular}
\caption{Convergence for the bean shape scatterer 
	 with $\kappa =2$ and $n = \sin(\pi x)\cos(\pi y) \hspace{1mm}\text{for} \hspace{1mm} \bm{x} = (x,y) \in \Omega$,
	 when $3$-point Newton-Cotes quadrature
	 employed for boundary patch integration in transverse integration.} \label{table:-Bean-K2-Q3}
\end{center}
\begin{center}
 \begin{tabular}{c|c|c|c|c|c|c}\hline\label{2}
Grid Size &  Unknown & Iteration & \multicolumn{2}{c|}{$L^2$} & \multicolumn{2}{c}{ $L^{\infty}$ } \\
 \cline{4-7}
 &   & Number & $\varepsilon_{2}$ & Order  & $\varepsilon_{\infty}$ &  Order \\ 
\hline
$2\times5\times9+ 1\times9\times9$ & 171 &3 & 5.25e-02 & & 1.69e-01 & - \\
\hline
$2\times9\times17+ 1\times17\times 17$ & 595 & 4 & 9.12e-03 & 2.53e+00 &3.34e-02 & 2.34e+00 \\
\hline
$2\times17\times33+ 1\times33\times 33$ & 2211 & 5 & 7.99e-04 & 3.51e+00 & 3.14e-03 & 3.41e+00 \\
\hline
$2\times33\times65+ 1\times65\times65$ & 8515 & 7& 4.51e-05 & 4.15e+00 & 8.14e-05 & 5.27e+00\\
\hline
$2\times65\times129+ 1\times129\times129$ & 33411 & 8 & 6.14e-07 & 6.20e+00 & 1.23e-06 & 6.04e+00\\
\hline
\end{tabular}
\caption{Convergence for the bean shape scatterer with $\kappa =2$ and $n = \sin(\pi x)\cos(\pi y) \hspace{1mm} \text{for} \hspace{1mm} \bm{x} = (x,y) \in \Omega$,
	 when $5$-point Newton-Cotes quadrature
	 employed for boundary patch integration in transverse integration.}\label{table:-Bean-K2-Q5}
\end{center}
\end{table} 

For the sake of pictorial visualization, results of scattering computation for $\kappa a  =20 $, $n =\sqrt{2}$
are shown in Figure (\ref{fig:-Bean-K10}) for a plane wave incidence traveling along positive $x$-axis. The plots visualize a solution obtained at the end of 82 iteration of GMRES when residual had reached $10^{-5}$. 
\end{exmp}

\section{An extension to three dimensions}
\label{sec::3d}
Though this text is primarily dedicated to presenting a fast and accurate computational strategy to solve the inhomogeneous acoustic scattering problem in two dimensions, as we mentioned in Section 2, this methodology has a straightforward extension
that allows for numerical solution of the corresponding three dimensional counterpart. The difficulty in high-order evaluation of the integral operator 
%(\ref{eq:-VolInt-B-I})
in three dimensions are largely analogous to those that appear in the two dimensional setting. Given that the main algorithmic steps remain unchanged, in this section, we avoid much of the repetitions and only briefly highlight some of the salient points underlying the extension. We then present a numerical verification of the fact that this methodology, indeed, produce rapidly converging numerical solutions to the inhomogeneous scattering problems.

\begin{table}[h!]
\begin{center}
\begin{tabular}{c|c|c|c|c|c|c}\hline\label{2}
Grid Size &  Unknown & Iter & \multicolumn{2}{c|}{$L^2$} & \multicolumn{2}{c}{ $L^{\infty}$ } \\
 \cline{4-7}
 &   &  & $\varepsilon_{2}$ & Order  & $\varepsilon_{\infty}$ &  Order \\ 
\hline
$2\times5\times5\times5+ 1\times5\times5\times5 $& 375 &  2& 2.03e-01 & - & 2.66e-01 & -  \\
\hline
$2\times9\times9\times9+ 1\times9\times9\times9 $& 2187 &  2&5.21e-02 & 1.96e+00&2.48e-01 & 1.01e-01 \\
\hline
$2\times17\times17\times17+ 1\times17\times17\times17 $& 14739 &  4&5.02e-03 & 3.38e+00&3.28e-02 & 2.92e+00 \\
\hline
$2\times33\times33\times33+ 1\times33\times33\times33 $& 107811 &  6&1.65e-04 & 4.93e+00&7.59e-04 & 5.44e+00 \\
\hline
$2\times65\times65\times65+ 1\times65\times65\times65 $& 823875 &  8 &4.55e-06 & 5.18e+00&1.61e-05 & 5.55e+00 \\
\hline
\end{tabular} \caption{Convergence study for scattering of plane wave  by a penetrable spherical shape scatterer 
	 with wave number $\kappa =2$ and refractive index $n = \sqrt{2}$, when $5$-point Newton-Cotes quadrature
	 employed for boundary patch integration in transverse integration.}\label{table3D:-K2}
\end{center}
\end{table}

We begin by recalling that, as in (\ref{eq:-VolInt-B-I}), %rewritten like equation (\ref{eq:-VolInt-B-II})
the integral $\mathcal{K}[\mathfrak{u}](\bm{x})$ is written as a sum of integrals over interior and boundary patches. In three dimensions, each of these integrals take the following form when expressed in the parametric space variables:
\[
\int \limits_{\mathcal{P}_{k}}
 G_{\kappa}(\bm{x},\bm{y}) m(\bm{y})\mathfrak{u}(\bm{y})\omega_{k}(\bm{y})\,d\bm{y} = 
 \iiint \limits_{[0,1]^3}
 G_{\kappa}(\bm{x},\bm{\xi}_k(t_1,t_2,t_3)) \varphi_k[\mathfrak{u}](t_1,t_2,t_3) \xi'_{k}(t_1,t_2,t_3) \,dt_1 \, dt_2 \, dt_3.
\]
As in the two dimensions, for the cases $k \not\in \mathcal{M}_B(\bm{x})$ and $k \not\in \mathcal{M}_I(\bm{x})$,  
the kernel $G_{\kappa}(\bm{x},\bm{y})$ remains non-singular within the integration region and, therefore, corresponding integrals can be approximated using high-order quadratures. 
%As we pointed out in two dimensional case,
%(see section 3.) it is advantageous to use different quadrature for both of the case. 
Again, when $k \not\in \mathcal{M}_I(\bm{x})$, the integrands have smooth and periodic extension
to $\mathbb{R}^{3}$ and the trapezoidal rule yields approximations with super-algebraic convergence. For the case when $k \not\in \mathcal{M}_B(\bm{x})$, on the other hand, high order can be attained by simply employing spectrally accurate 
Trapezoidal rule for planar integrations (with respect to $t_1, t_2$) and a high-order composite Newton-Cotes quadrature for the integration in the transverse variable $t_3$.

For the cases $k \in \mathcal{M}_B(\bm{x})$ and $k \in \mathcal{M}_I(\bm{x})$, as before, high-order
evaluation of the integrals becomes difficult owing to the kernel singularity at $\bm{x} = \bm{y}$.  
%Again, similar to the two dimensional
%case, our strategy for dealing with these integrand singularity differ for boundary and interior patches. 
Toward dealing with this in boundary patches, following \cite{anand2007efficient}, we change to polar coordinates $(\rho,\theta)$ centered around $\bar{t}^{\bm{x}} = (t_{1}^{\bm{x}},t_{2}^{\bm{x}})$, the projection of $\xi_{k}^{-1}(\bm{x})$ on
to the $t_3$-integration plane, %$(\bar{t}^{\bm{x}},t_{3}^{\bm{x}})$, 
which, upon employing an accompanying polynomial change of variable in $\rho$, provides an effective resolution of difficulties arising out of the kernel singularity. 
Subsequent application of Trapezoidal rule in $t_1$, $t_2$ variables and a composite Newton-Cotes quadrature for $t_3$ integration yields accurate results.
%We refer the readers to \cite{anand2007efficient} for a more detailed discussion on this approach. 
For the singular integral over interior patches, on the other hand, we change to spherical coordinates $(\rho,\theta,\phi)$ around target point $(t_1^{\bm{x}},t_2^{\bm{x}},t_3^{\bm{x}})$. This, in turn, yields smooth integrands, where accurate integrations can be affected by employing Trapezoidal rule in  $\rho,\phi$ variables and Clenshaw-€"Curtis quadrature in $\theta$ variable.

The integration scheme can again be accelerated by a suitable use of two face equivalent source approximations on Cartesian grids. Just as in two dimensions, this strategy of employing three dimensional FFTs for approximating the convolution, further reduces the evaluation cost of the 
non-singular non-adjacent interactions. In this case, the algorithm exhibits the theoretical computational complexity of $O(N \log N)$ with respect to the grid size $N$, provided we choose $L = O(N^{1/3})$, when $L^3$ number of sub-cubes are used to cover the inhomogeneity.

\begin{figure}[h!]
\begin{center}
\subfigure[ Spherical shape scatterer]{\includegraphics[ clip=true,trim=500 50  450 45, scale =.36]{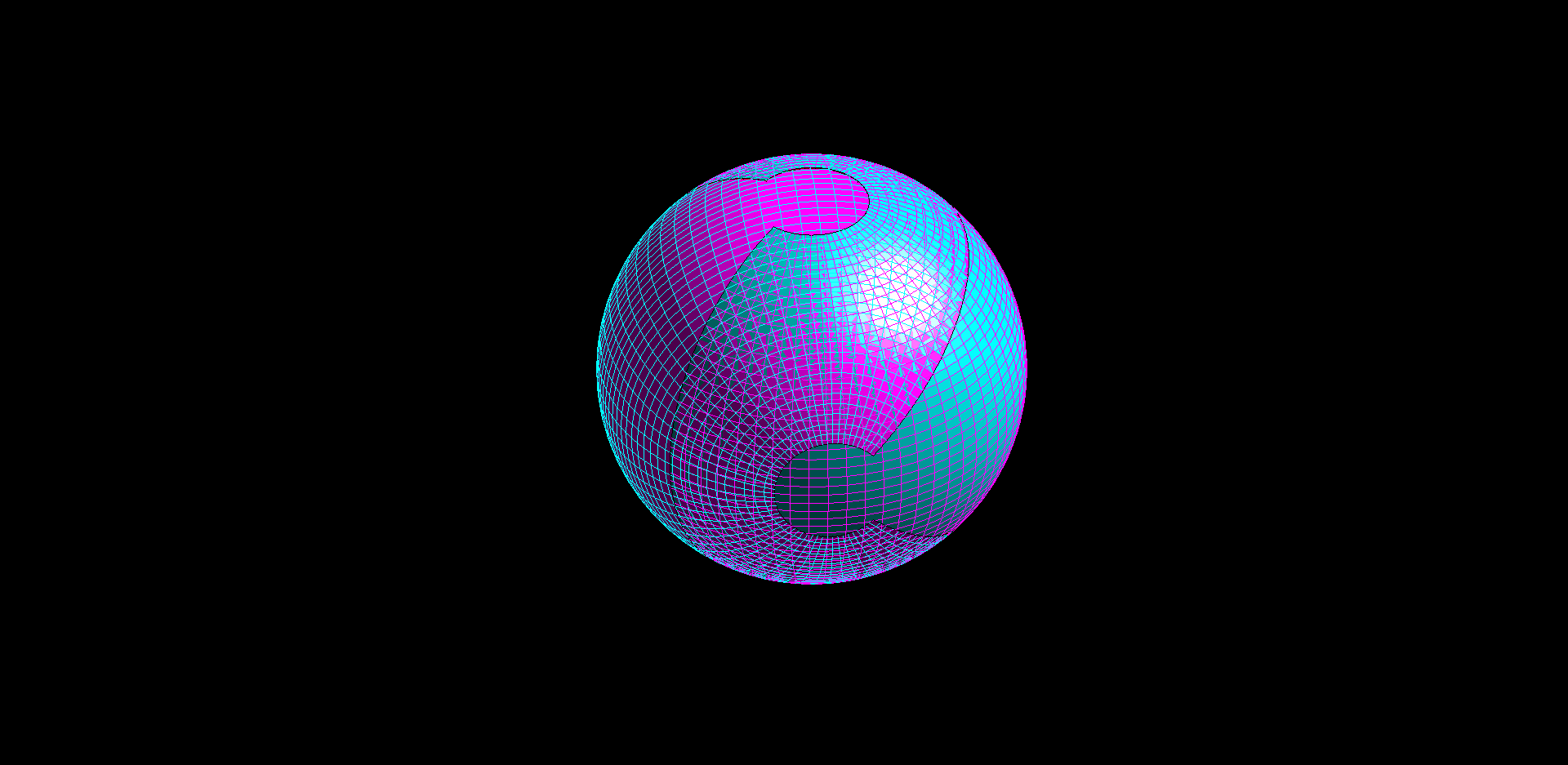}}
\hfill
\subfigure[Intensity of the total field ($|\mathfrak{u}|^2$)]{\includegraphics[clip=true, trim=450 50  400 45, scale=.36]{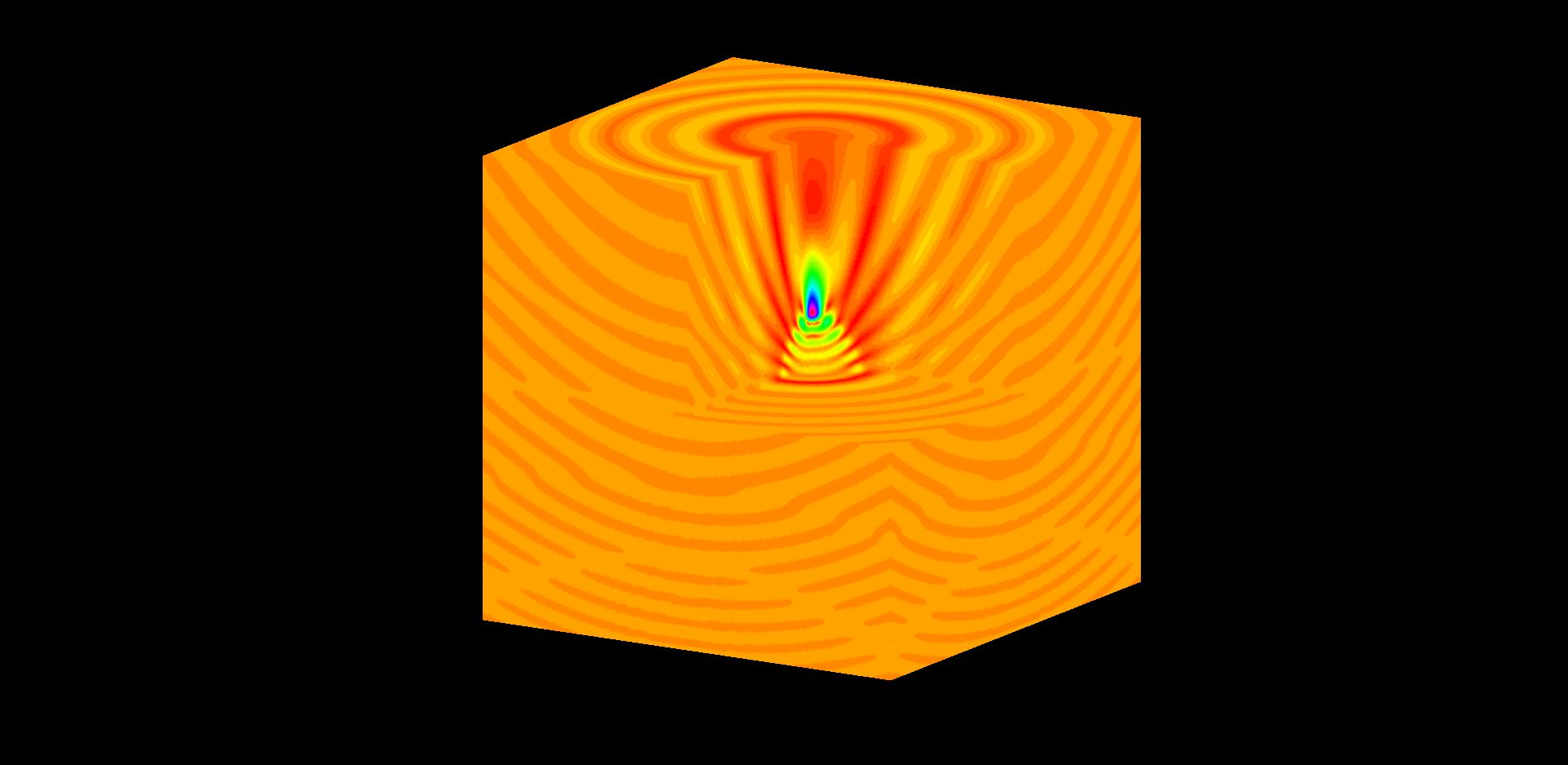}}
\caption{The figure of the left depicts a the discretization on the outer surface of the scatterer coming from a typical volumetric computational grid. The figure of the right presents a visualization of the total field intensity.} \label{fig:-Sphere-K5}
\end{center}
\end{figure}

%\subsection{Numerical Results}
In order to demonstrate the high-order convergence of the three dimensional scheme, we present %some numerical results for 
results from a plane wave scattering by a penetrable spherical inhomogeneity with constant refractive index, $n =\sqrt{2}$. 
%Since for spherical scatterer exact solution is known therefore 
%computed solution have been compared against the exact solution. We assume that incidence plane wave impinges on the obstacle from 
%positive $z$-axis.  
The convergence study, given in Table \ref{table3D:-K2}, that corresponds to the incident plane wave $\mathfrak{u}^i = \exp(i\kappa z)$ with wavenumber $\kappa = 2$, clearly exhibits rapid convergence of numerical solutions.

% \begin{table}[h!]
% \begin{center}
% \begin{tabular}{c|c|c|c|c|c|c}\hline\label{2}
% Grid Size &  Unknown & Iter &  $\varepsilon_{2}$ & Order  & $\varepsilon_{\infty}$ &  Order \\ 
% \hline
% $2\times5\times5\times5+ 1\times5\times5\times5 $& 375 &  2& 2.03e-01 & - & 2.66e-01 & -  \\
% \hline
% $2\times9\times9\times9+ 1\times9\times9\times9 $& 2187 &  2&5.21e-02 & 1.96e+00&2.48e-01 & 1.01e-01 \\
% \hline
% $2\times17\times17\times17+ 1\times17\times17\times17 $& 14739 &  4&5.02e-03 & 3.38e+00&3.28e-02 & 2.92e+00 \\
% \hline
% $2\times33\times33\times33+ 1\times33\times33\times33 $& 107811 &  6&1.65e-04 & 4.93e+00&7.59e-04 & 5.44e+00 \\
% \hline
% $2\times65\times65\times65+ 1\times65\times65\times65 $& 823875 &  8 &4.55e-06 & 5.18e+00&1.61e-05 & 5.55e+00 \\
% \hline
% \end{tabular} \caption{Convergence study for scattering of plane wave  by a penetrable spherical shape scatterer 
% 	 with wave number $\kappa =2$ and refractive index $n = \sqrt{2}$, when $5$-point Newton-Cotes quadrature
% 	 employed for boundary patch integration in transverse integration.}\label{table3D:-K2}
% \end{center}
% \end{table}

As a final example, we visualize, in Figure \ref{fig:-Sphere-K5}, the results of a scattering computation corresponding to a spherical obstacle of acoustical size $ \kappa a =10$ and refractive index $n=2$.
For this simulations, a relative error of order $10^{-5}$ is achieved at a computational grid of size $2\times 33 \times 65 \times 65+1\times 65 \times 65 \times 65 $.

\section{Conclusions}
\label{sec::conclu}
In this paper, we discuss a fast high-order method for scattering of acoustic waves by penetrable
inhomogeneous media which require $O(N\log N)$ computational cost for each iteration of iterative linear system solvers. 
Rapidly convergent approximations are obtained through a use of specialized quadrature rule while a reduced computational cost is achieved by utilizing a strategy based on equivalent sources
approximation on the Cartesian grids. We present a series of numerical experiments to demonstrate its performance in terms of computational efficiency as well as numerical accuracy.
We emphasize that this algorithm is designed to retain high-order convergence even when material properties
jump across the scattering interface and can be employed for simulating scattering of acoustic waves by geometrically complex structures. 
For instance, in our numerical experiments where function $m =1-n^{2}$ is discontinuous across
the material interface which, of course, limits the global smoothness of the total field $\mathfrak{u}$, (in fact, $\mathfrak{u}\in C^{1}(\mathbb{R}^{3})$ \cite{potthast2001point}), we still obtain rapidly convergent approximations. This, of course, becomes possible as a result of carefully avoiding integrating across the material interface.

%One can see in numerical results that
%the number of GMRES iteration increases quadratically  as the problem size increase for instance, {\color{red} in example 4 (section 3.3) 
%for complex bean shape scatterer of acoustical size $\kappa a = 20$, 
%$82$ number of GMRES iterations is needed to solve the linear system.} Although in recent years, important efforts have been directed toward
%the design of preconditioning techniques \cite{bruno2004efficient,sifuentes2010preconditioned,ying2014sparsifying} to reduce number of iteration but problem still can not be considered completely solved and we are
%intend to develop efficient preconditioners for the integral equation formulation. Another future work along this direction is
%to apply our method for electromagnetic and elastic scattering problems. Further we wish to generalize our method for non-smooth
%geometry such as scatterer containing corner, edge and cusps.

\section*{Acknowledgements}

Akash Anand gratefully acknowledges support from SERB-DST through contract No.  SERB/F/ \\5152/2013-2014. Ambuj Pandey gratefully acknowledges support from CSIR.



%\section*{References}

\bibliography{LS2D.bib}
\bibliographystyle{abbrv}
%\bibliographystyle{elsarticle-num} 


%\begin{thebibliography}{00}


%\bibitem{}

%\end{thebibliography}
\end{document}